\documentclass{article}
\usepackage{amssymb}
\usepackage{amsmath}
\usepackage{tikz}
\usepackage{graphicx}
\usepackage{float}

\newtheorem{theorem}{Theorem}

\newtheorem{lemma}[theorem]{Lemma}

\newtheorem{remark}[theorem]{Remark}

\newcommand{\func}[1]{\text{#1}}
\begin{document}

\title{\textbf{{\Large Backward Uniqueness for a PDE Fluid-Structure
Interaction{\ }}}}
\author{ George Avalos \\
%EndAName
{\normalsize {University of Nebraska-Lincoln}}\\
{\small {Lincoln NE, U.S.A.}}\\
{\small \texttt{gavalos@math.unl.edu}} \and Thomas Clark \\
%EndAName
{\normalsize {Dordt College}}\\
{\small {Sioux Center IA, U.S.A.}}\\
{\small \texttt{Tom.Clark@dordt.edu}}}
\date{July 21, 2014}
\maketitle

\begin{abstract}
In this work, we establsh the so-called backward unqiueness property for a
coupled system of partial differential equations (PDEs) which governs a
certain fluid-structure interaction. In particular, a three-dimensional
Stokes flow interacts across a boundary interface with a two-dimensional
mechanical plate equation. By way of attaining this result, a certain
estimate is obtained for the associated semigroup generator resolvent.
\end{abstract}

%\date{February 25, 2014}

%%%%%%%%%%%%%%%%%%%%%%%%%%%%%%
%
%
%In \cite{AC} it was shown that the (\ref{1})-(\ref{ic}) is equivalent to solving the
%operator equation
%\begin{equation}
%\begin{bmatrix} w(t) \\ w_t(t) \\ u(t) \end{bmatrix} = e^{\mathcal{A}_\rho t}
%\begin{bmatrix} w_1^{\ast} \\ w_2^{\ast} \\u_0^{\ast} \end{bmatrix}
%\end{equation}
%for a particular operator $\mathcal{A_\rho} : D(\mathcal{A_\rho} \subset \mathbf{H}_\rho
%\rightarrow \mathbf{H}_\rho$.  The pressure $p(t)$ is given pointwise in time via
%\begin{equation}
%p(t) = G_{\rho, 1} (w(t)) + G_{\rho, 2}(u(t));
%\end{equation}
%see \cite{AC} for specific details about these operators.
%
%
%%%%%%%%%%%%%%%%%%%%%%%%%%

\section{Statement of the Problem and Main Result}

We consider here the problem of establishing the so-called backward
uniqueness property for the partial differential equation (PDE) model given
in \cite{igor} and \cite{AC}, which describes a certain fluid-structure
interactive dynamics. One novelty of this PDE system is the unique way in
which the geometry affects the coupling between the fluid and the plate.
Since the coupling involves the pressure term, the system cannot be solved
via the classic Leray projector. Instead in \cite{AC} wellposedness is given
via a semigroup formulation and proved via the Babu\v{s}ka-Brezzi theorem.

As explained in more detail in \cite{AT}, the backward uniqueness property
-- described below in Theorem \ref{back} -- has important implications for
the controllability of the system in the sense of PDE control theory.

As mentioned above, the geometry plays are particular role in the
wellposedness of the system. The fluid chamber $\mathcal{O}\subset \mathbb{R}%
^{3}$ will be a bounded domain with sufficiently smooth boundary. Moreover, $%
\partial \mathcal{O}=\bar{\Omega}\cup \bar{S}$, with $\Omega \cap
S=\varnothing $. More specifically%
\begin{equation*}
\Omega \subset \left\{ x=(x_{1,}x_{2},0)\right\} \text{, and surface }%
S\subset \left\{ x=(x_{1,}x_{2},x_{3}):x_{3}\leq 0\right\} .
\end{equation*}%
In consequence, if $\nu (x)$ denotes the exterior unit normal vector to $%
\partial \mathcal{O}$, then%
\begin{equation}
\left. \nu \right\vert _{\Omega }=[0,0,1]  \label{normal}
\end{equation}%

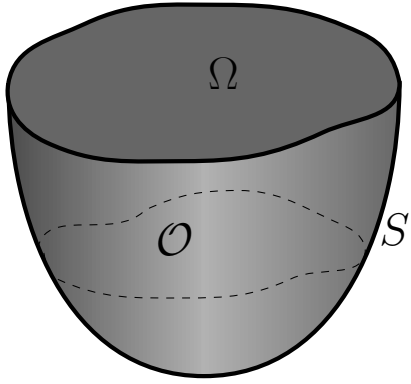
\begin{figure}[H]
\begin{center}
\begin{tikzpicture}[scale=1.3]
\draw[left color=black!70,right color=black!60,middle color=black!30, ultra thick] (-2,0,0) to [out=0, in=180] (2,0,0) to [out=270, in = 0]
(0,-3,0) to [out=180, in =270] (-2,0,0);

\draw [fill=black!60, ultra thick] (-2,0,0) to [out=80, in=205](-1.214,.607,0) to [out=25, in=180](0.2,.8,0) to [out=0, in=155] (1.614,.507,0)
to [out=335, in=100](2,0,0) to [out=270, in=25] (1.214,-.507,0) to [out=205, in=0](-0.2,-.8,0) [out=180, in=335] to (-1.614,-.607,0) to
[out=155, in=260] (-2,0,0);

\draw [dashed, thin] (-1.7,-1.7,0) to [out=80, in=225](-.6,-1.3,0) to [out=25, in=180](0.35,-1.1,0) to [out=0, in=155] (1.3,-1.4,0) to [out=335,
in=100](1.65,-1.7,0) to [out=270, in=25] (0.9,-2.0,0) to [out=205, in=0](-0.2,-2.2,0) [out=180, in=335] to (-1.514,-2.0) to [out=155, in=290]
(-1.65,-1.7,0);

\node at (0.2,0.1,0) {{\LARGE$\Omega$}};

\node at (1.95,-1.5,0) {{\LARGE $S$}};

\node at (-0.3,-1.6,0) {{\LARGE $\mathcal{O}$}};
\end{tikzpicture}
\caption{The Fluid-Structure Geometry}
\end{center}
\end{figure}

In addition, $\left[ \mathcal{O},\Omega \right] $ is assumed to fall within
one of the following classes:%
\begin{equation*}
\begin{array}{l}
\mathsf{(G.1)}\ \mathcal{O}\text{ is a convex domain with wedge angles }<%
\frac{2\pi }{3}\text{. \ Moreover, }\Omega \\
\text{ \ \ \ \ \ \ \ has smooth boundary, and }S\text{ is a piecewise smooth
surface;} \\
\mathsf{(G.2)}\ \mathcal{O}\text{ is a convex polyhedron having angles}<%
\frac{2\pi }{3}\text{,} \\
\text{ \ \ \ \ \ \ \ \ \ \ \ \ \ \ and so then }\Omega \text{ is a convex
polygon with angles}<\frac{2\pi }{3}\text{.}%
\end{array}%
\end{equation*}

The PDE model is as follows, with \textquotedblleft rotational inertia
parameter\textquotedblright\ $\rho \geq 0$, and in solution variables $%
[w(x,t),w_{t}(x,t)]$, $u(x,t)=[u^{1}(x,t),u^{2}(x,t),u^{3}(x,t)]$, and $%
p(x,t)$:
\begin{align}
& w_{tt}-\rho \Delta w_{tt}+\Delta ^{2}w=\left. p\right\vert _{\Omega }\text{
\ in }\Omega \times (0,T),  \label{1} \\
& w=\frac{\partial w}{\partial \nu }=0\text{ \ on }\partial \Omega ;
\label{2} \\
& u_{t}-\Delta u+\nabla p=0\text{ \ in }\mathcal{O}\times (0,T),  \label{3}
\\
& \mathrm{div}(u)=0\text{ \ in }\mathcal{O}\times (0,T),  \label{4} \\
& u=0\text{ on }S\text{ \ and }u=[u^{1},u^{2},u^{3}]=[0,0,w_{t}]\text{ \ on }%
\Omega ,  \label{5}
\end{align}%
with initial conditions%
\begin{equation}
\lbrack w(0),w_{t}(0),u(0)]=[w_{1}^{\ast },w_{2}^{\ast },u_{0}^{\ast }]\in
\mathbf{H}_{\rho }\text{.}  \label{ic}
\end{equation}%
Here, the finite energy space $\mathbf{H}_{\rho }$ is given by
\begin{eqnarray}
\mathbf{H}_{\rho } &=&\Big\{\left[ \omega _{1},\omega _{2},f\right] \in %
\left[ H_{0}^{2}(\Omega )\cap \widehat{L}^{2}(\Omega )\right] \times W_{\rho
}\times \mathcal{H}_{\mathrm{fluid}}  \notag \\
&&\text{ \ \ \ \ }\text{with }\left. f\cdot \nu \right\vert _{\Omega
}=[0,0,f^{3}]\cdot \lbrack 0,0,1]=\omega _{2}\Big\},  \label{energy}
\end{eqnarray}%
where%
\begin{equation}
\widehat{L}^{2}(\Omega )=\left\{ \varpi \in L^{2}(\Omega ):\int_{\Omega
}\varpi d\Omega =0\right\} ;  \label{L_hat}
\end{equation}%
\begin{equation}
\mathcal{H}_{\mathrm{fluid}}=\left\{ f\in \mathbf{L}^{2}(\mathcal{O}):%
\mathrm{div}(f)=0\text{; }\left. f\cdot \nu \right\vert _{S}=0\right\} ;
\label{H_f}
\end{equation}%
and%
\begin{equation}
W_{\rho }=\left\{
\begin{array}{l}
\widehat{L}^{2}(\Omega )\text{, \ if }\rho =0, \\
\\
H_{0}^{1}(\Omega )\cap \widehat{L}^{2}(\Omega )\text{, \ if }\rho >0.%
\end{array}%
\right.  \label{velocity}
\end{equation}%
This Hilbert space $\mathbf{H}_{\rho }$ of finite energy is endowed here
with the followng norm-inducing inner product:%
\begin{equation}
\left( \big[\omega _{1},\omega _{2},f\big],\big[\tilde{\omega}_{1},\tilde{%
\omega}_{2},\tilde{f}\big]\right) _{\mathbf{H}_{\rho }}=(\Delta \omega
_{1},\Delta \tilde{\omega}_{1})_{\Omega }+(\omega _{2},\tilde{\omega}%
_{2})_{\Omega }+\rho (\nabla \omega _{2},\nabla \tilde{\omega}_{2})_{\Omega
}+(f,\tilde{f})_{\mathcal{O}}.  \label{norm}
\end{equation}%
where $(\cdot ,\cdot )_{\Omega }$ and $(\cdot ,\cdot )_{\mathcal{O}}$ are
the $L^{2}$-inner products on their respective geometries.

\medskip

For the PDE system (\ref{1})-(\ref{ic}), semigroup well-posedness result was
established in \cite{AC}; a proof of wellposedness, for $\rho >0$, via a
Galerkin method was also given in \cite{igor}, which paper was primarily
concerned with longtime behaviour of corresponding solutions, under
nonlinear effects. We will presently give an explicit description of the
modeling semigroup generator $\mathcal{A}_{\rho }:D(\mathcal{A}_{\rho
})\subset \mathbf{H}_{\rho }\rightarrow \mathbf{H}_{\rho }$; its
construction in \cite{AC} hinged upon an appropriate elimination of the
pressure variable $p$ in (\ref{1})-(\ref{ic}). (As the no-slip boundary
condition is necessarily not in play for the fluid variable $u$, one cannot
merely invoke the classic Leray projector to eliminate the pressure term, as
one would in uncoupled fluid flow; see e.g., \cite{temam}.)

\begin{theorem}
\label{well}(See \cite{AC}.) The PDE model (\ref{1})-(\ref{ic}) is
associated with a $C_{0}$- contraction semigroup $\left\{ e^{\mathcal{A}%
_{\rho }t}\right\} _{t\geq 0}\subset \mathcal{L}(\mathbf{H_{\rho }}$), the
generator of which is given below in (\ref{A})-(\ref{domain}). Therewith,
for any initial data $[w_{1}^{\ast },w_{2}^{\ast },u_{0}^{\ast }]\in \mathbf{%
H_{\rho }}$, the solution $[w,w_{t},u]\in C([0,T;\mathbf{H_{\rho }})$ is
given by
\begin{equation*}
\left[
\begin{array}{c}
w(t) \\
w_{t}(t) \\
u(t)%
\end{array}%
\right] =e^{\mathcal{A}_{\rho }t}\left[
\begin{array}{c}
w_{1}^{\ast } \\
w_{2}^{\ast } \\
u^{\ast }%
\end{array}%
\right] .
\end{equation*}
\end{theorem}

\bigskip

The main result of this paper -- Theorem \ref{back} below -- deals with
establishing the aforesaid backward uniqueness property for the contraction $%
C_{0}$-semigroup associated with the PDE model (\ref{1})-(\ref{ic}). The
driving agent of our proof of Theorem \ref{back} is the following abstract
resolvent criterion for backward uniqueness.

\bigskip

\begin{theorem}
\label{semi}(See \cite{LRT}, Theorem 3.1, p. 225.) Let $A$ be the
infinitesimal generator of a s.c. semigroup $e^{At}$ in a Banach space $X$.
Assume that there exist constants $a\in (\pi /2,\pi )$, $r_{0}>0$, and $C>0$%
, such that%
\begin{equation*}
\left\Vert \mathcal{R}(re^{\pm ia};A)\right\Vert _{\mathcal{L}%
(X)}=\left\Vert (re^{\pm ia}I-A)^{-1}\right\Vert _{\mathcal{L}(X)}\leq C%
\text{,}
\end{equation*}%
for all $r\geq r_{0}$. Then the backward uniqueness property holds true;
that is, $e^{AT}x_{0}=0$ for $T>0$, $x_{0}\in X$, implies $x_{0}=0$.
\end{theorem}

\bigskip

By way of applying the abstract Theorem \ref{semi} to the modeling generator
$\mathcal{A}_{\rho }:D(\mathcal{A}_{\rho })\subset \mathbf{H}_{\rho
}\rightarrow \mathbf{H}_{\rho }$ of (\ref{1})-(\ref{ic}), given explicitly
in (\ref{A})-(\ref{domain}) below, we will consider the following resolvent
relation with complex parameter $\lambda =\alpha +i\beta $, which is
formally a \textquotedblleft frequency domain\textquotedblright\ version of (%
\ref{1})-(\ref{ic}):
\begin{equation}
\left( \lambda I-\mathcal{A}_{\rho }\right) \left[
\begin{array}{c}
\omega _{1} \\
\omega _{2} \\
\mu%
\end{array}%
\right] =\left[
\begin{array}{c}
\omega _{1}^{\ast } \\
\omega _{2}^{\ast } \\
\mu ^{\ast }%
\end{array}%
\right] \in \mathbf{H}_{\rho }.  \label{abstract}
\end{equation}%
Here, the pre-image $\left[ \omega _{1},\omega _{2},\mu \right] \in D(%
\mathcal{A}_{\rho })$ and forcing term $\left[ \omega _{1}^{\ast },\omega
_{2}^{\ast },\mu ^{\ast }\right] \in \mathbf{H}_{\rho }$.

\bigskip

With respect to the frequency domain parameter, we will furthermore impose
that $\lambda =\alpha +i\beta $ should obey the following criteria:

\begin{enumerate}
\item[Criterion 1:] $\lambda =\alpha +i\beta =\left\vert \lambda \right\vert
e^{\pm i\vartheta }$, for fixed $\vartheta \in \left( \frac{3\pi }{4},\pi
\right) \,$. (And so on either of these two rays, we have $0<\left\vert \tan
\vartheta \right\vert <1$, $\left\vert \beta \right\vert =\left\vert \alpha
\right\vert \left\vert \tan \vartheta \right\vert $, $\left\vert \lambda
\right\vert ^{2}=\alpha ^{2}+\beta ^{2}=\alpha ^{2}(1+\tan ^{2}\vartheta )$.)

\item[Criterion 2:] $\left\vert \alpha \right\vert >0$ is sufficiently large.
\end{enumerate}

\bigskip

Our main result can now be stated as follows:

\begin{theorem}
\label{back}(i) With respect to the resolvent relation (\ref{abstract}), or
the equivalent fluid-structure PDE (\ref{s})-(\ref{f}) below, let the
Criteria 1 and 2 be in force. Then for all $\rho \geq 0$, the solution $%
\left[ \omega _{1},\omega _{2},\mu \right] \in D(\mathcal{A}_{\rho })$ obeys
the following bound, which is uniform for all $\lambda =\alpha +i\beta
=\left\vert \lambda \right\vert e^{\pm i\vartheta }$, with fixed $\vartheta
\in \left( \frac{3\pi }{4},\pi \right) $, and $\left\vert \alpha \right\vert
>0$ sufficiently large:
\begin{equation}
\left\Vert \left[
\begin{array}{c}
\omega _{1} \\
\omega _{2} \\
\mu%
\end{array}%
\right] \right\Vert _{\mathbf{H}_{\rho }}\leq C_{\vartheta }\left\Vert \left[
\begin{array}{c}
\omega _{1}^{\ast } \\
\omega _{2}^{\ast } \\
\mu ^{\ast }%
\end{array}%
\right] \right\Vert _{\mathbf{H}_{\rho }}.  \label{bound}
\end{equation}%
(ii) In consequence, this estimate and Theorem \ref{semi} yields the
conclusion that the fluid-structure $C_{0}$-contraction semigroup $\left\{
e^{\mathcal{A}_{\rho }t}\right\} _{t\geq 0}$ satisfies the backward
uniqueness property: Namely, if for given $T>0$ and $[w_{1}^{\ast
},w_{2}^{\ast },u_{0}^{\ast }]\in \mathbf{H}_{\rho }$, one has
\begin{equation*}
e^{\mathcal{A}_{\rho }T}\left[
\begin{array}{c}
w_{1}^{\ast } \\
w_{2}^{\ast } \\
u_{0}^{\ast }%
\end{array}%
\right] =\vec{0}\text{, then necessarily }\left[
\begin{array}{c}
w_{1}^{\ast } \\
w_{2}^{\ast } \\
u_{0}^{\ast }%
\end{array}%
\right] =\vec{0}.
\end{equation*}
\end{theorem}

\bigskip

\begin{remark}
Unlike the coupled PDE examples in \cite{LRT}, \cite{AT}, \cite{AT2}, the
frequency domain estimate (\ref{bound}), for the fluid-structure solution $%
\left[ \omega _{1},\omega _{2},\mu \right] $ of (\ref{abstract}), does not
manifest a uniform rate of decay with respect to $\func{Re}\lambda =\alpha $%
. We are not certain that such a decay is actually possible. However by
Theorem \ref{semi}, the uniform bound (\ref{bound}) suffices to established
the aforesaid backward uniquess property.
\end{remark}

\begin{remark}
In the course of proof, the reader could infer that for the rotational
inertial case $\rho >0$, one will in fact have the uniform estimate (\ref%
{bound}) for any rays along the angle $\vartheta \in \left( \frac{\pi }{2}%
,\pi \right) $, $\left\vert \lambda \right\vert $ large enough.
\end{remark}

\bigskip

\section{The Description of the Fluid-Structure Generator}

Under the geometric conditions \textsf{(G.1)} and \textsf{(G.2)}, we now
tersely define the modeling generator $\mathcal{A}_{\rho }:D(\mathcal{A}%
_{\rho })\subset \mathbf{H}_{\rho }\rightarrow \mathbf{H}_{\rho }$ which
describes the dynamics (\ref{1})-(\ref{ic}), and for which Theorem \ref{back}
applies. Full details are given in \cite{AC} and \cite{AB}.

\medskip

To start, let $A_{D}:L^{2}(\Omega )\rightarrow L^{2}(\Omega )$ be given by%
\begin{equation}
A_{D}g=-\Delta g\text{, \ \ }D(A_{D})=H^{2}(\Omega )\cap H_{0}^{1}(\Omega ).
\label{dirichlet}
\end{equation}%
If we subsequently make the denotation for all $\rho \geq 0$,%
\begin{equation}
P_{\rho }=I+\rho A_{D}\text{, \ }D(P_{\rho })=\left\{
\begin{array}{l}
L^{2}(\Omega )\text{, \ if }\rho =0, \\
D(A_{D})\text{, \ if }\rho >0,%
\end{array}%
\right.  \label{P}
\end{equation}%
then the mechanical PDE component (\ref{1})-(\ref{2}) can be written as
\begin{equation*}
P_{\rho }w_{tt}+\Delta ^{2}w=\left. p\right\vert _{\Omega }\text{ on\ }(0,T).
\end{equation*}%
Using the characterization from \cite{grisvard} that%
\begin{equation*}
\text{\ }D(P_{\rho }^{\frac{1}{2}})=\left\{
\begin{array}{l}
L^{2}(\Omega )\text{, \ if }\rho =0, \\
H_{0}^{1}(\Omega )\text{, \ if }\rho >0,%
\end{array}%
\right.
\end{equation*}%
then from (\ref{norm}) we can rewrite%
\begin{equation}
\left( \left[ \omega _{1},\omega _{2},f\right] ,\left[ \tilde{\omega}_{1},%
\tilde{\omega}_{2},\tilde{f}\right] \right) _{\mathbf{H}_{\rho }}=(\Delta
\omega _{1},\Delta \tilde{\omega}_{1})_{\Omega }+(P_{\rho }^{\frac{1}{2}%
}\omega _{2},P_{\rho }^{\frac{1}{2}}\tilde{\omega}_{2})_{\Omega }+(f,\tilde{f%
})_{\mathcal{O}}.  \label{norm1}
\end{equation}

\medskip

Moreover, in order to eliminate the pressure -- see \cite{AC} -- we require
the following \textquotedblleft Robin\textquotedblright\ maps $R_{\rho }$
and $\tilde{R}_{\rho }$:%
\begin{eqnarray}
R_{\rho }g &=&f\Leftrightarrow \left\{ \Delta f=0\text{ \ in }\mathcal{O};%
\text{ }\frac{\partial f}{\partial \nu }+P_{\rho }^{-1}f=g\text{ \ on }%
\Omega \text{; \ }\frac{\partial f}{\partial \nu }=0\text{ on }S\right\} .
\label{R1} \\
&&  \notag \\
\tilde{R}_{\rho }g &=&f\Leftrightarrow \left\{ \Delta f=0\text{ \ in }%
\mathcal{O};\text{ }\frac{\partial f}{\partial \nu }+P_{\rho }^{-1}f=0\text{
\ on }\Omega \text{; \ }\frac{\partial f}{\partial \nu }=g\text{ on }%
S\right\} .  \label{R2}
\end{eqnarray}%
By Lax-Milgram we then have%
\begin{equation}
R_{\rho }\in \mathcal{L}\big(H^{-\frac{1}{2}}(\Omega ),H^{1}(\mathcal{O})%
\big)\text{; \ }\tilde{R}_{\rho }\in \mathcal{L}\big(H^{-\frac{1}{2}%
}(S),H^{1}(\mathcal{O})\big).  \label{Rs}
\end{equation}%
(We are also using implicity the fact that $P_{\rho }^{-1}$ is positive
definite, self-adjoint on $\Omega $.)

Therewith, it is shown in \cite{AC} that the pressure variable $p(t)$ can be
written pointwise in time as
\begin{equation}
p(t)=G_{\rho ,1}(w(t))+G_{\rho ,2}(u(t)),  \label{p}
\end{equation}%
where
\begin{eqnarray}
G_{\rho ,1}(w) &=&R_{\rho }(P_{\rho }^{-1}\Delta ^{2}w);  \label{G1} \\
G_{\rho ,2}(u) &=&R_{\rho }(\left. \Delta u^{3}\right\vert _{\Omega })+%
\tilde{R}_{\rho }(\left. \Delta u\cdot \nu \right\vert _{S}).  \label{G2}
\end{eqnarray}

With these operators, we defined in \cite{AC} the generator $\mathcal{A}%
_{\rho }:D(\mathcal{A}_{\rho })\subset \mathbf{H}_{\rho }\rightarrow \mathbf{%
H}_{\rho }$, which is associated with the fluid structure system (\ref{1})-(%
\ref{ic}):
\begin{eqnarray}
&&\mathcal{A}_{\rho }\equiv
\begin{bmatrix}
0 & I & 0 \\
-P_{\rho }^{-1}\Delta ^{2}+P_{\rho }^{-1}G_{\rho ,1}\big|_{\Omega } & 0 &
P_{\rho }^{-1}G_{\rho ,2}\big|_{\Omega } \\
-\nabla G_{\rho ,1} & 0 & \Delta -\nabla G_{\rho ,2}%
\end{bmatrix}%
;  \label{A} \\
&&  \notag \\
&&\text{with }D(\mathcal{A}_{\rho })=\big\{\left[ w_{1},w_{2},u\right] \in
\mathbf{\ H}_{\rho }\text{ satisfying}:  \notag \\
&&\text{ \quad \quad (a) }w_{1}\in \mathcal{S}_{\rho }\equiv \left\{
\begin{array}{l}
H^{4}(\Omega )\cap H_{0}^{2}(\Omega )\text{, \ if }\rho =0; \\
H^{3}(\Omega )\cap H_{0}^{2}(\Omega )\text{, \ if }\rho >0;%
\end{array}%
\right.  \label{S} \\
&&\text{ \quad \quad (b) }w_{2}\in H_{0}^{2}(\Omega )\text{, }u\in \mathbf{H}%
^{2}(\mathcal{O});  \notag \\
&&\text{ \quad \quad (c) }u=\vec{0}\text{ on }S\text{ and }u=[0,0,w_{2}]%
\text{ on }\Omega \big\}.  \label{domain}
\end{eqnarray}

\bigskip

\begin{remark}
Given data $\left[ w_{1},w_{2},u\right] \in D(\mathcal{A}_{\rho })$, note
that as $\Delta u\in L^{2}(\mathcal{O})$ and $\text{div}(\Delta u)=0$, then
by Theorem 1.2, p. 9 in \cite{temam}, we have the trace regularity
\begin{equation}
\Delta u\cdot \nu \big|_{\partial \mathcal{O}}\in H^{-\frac{1}{2}}(\partial
\mathcal{O});  \label{CritReg}
\end{equation}%
and so the pressure term
\begin{equation}
p\equiv G_{\rho ,1}(w_{1})+G_{\rho ,2}(u)\in H^{1}(\mathcal{O}).
\label{pInH1}
\end{equation}%
Thus, $\mathcal{A}_{\rho }:D(\mathcal{A}_{\rho })\subset \mathbf{H}_{\rho
}\rightarrow \mathbf{H}_{\rho }$ is indeed well-defined (see in particular
the $2-3$ and $3-3$ entries of matrix $\mathcal{A}_{\rho }$).
\end{remark}

\bigskip

It is shown in \cite{AC} that $\mathcal{A}_{\rho }:D(\mathcal{A}_{\rho
})\subset \mathbf{H}_{\rho }\rightarrow \mathbf{H}_{\rho }$ is maximal
dissipative, thereby giving rise to Theorem \ref{well} above. The next
Section is devoted to the proof of the main, backward uniqueness result.

\medskip

\section{Proof of Theorem \protect\ref{back}}

\medskip

With $\lambda =\alpha +i\beta $, and with the definition of $\mathcal{A}%
_{\rho }:D(\mathcal{A}_{\rho })\subset \mathbf{H}_{\rho }\rightarrow \mathbf{%
H}_{\rho }$ in hand from (\ref{A})-(\ref{domain}), the resolvent relation (%
\ref{abstract}) gives rise to the following fluid-structure PDE system:%
\begin{align}
& \left\{
\begin{array}{l}
\omega _{2}=\lambda \omega _{1}-\omega _{1}^{\ast }\text{ \ in }\Omega \\
\\
\left( \alpha ^{2}-\beta ^{2}\right) \omega _{1}+2i\alpha \beta \omega
_{1}+P_{\rho }^{-1}\Delta ^{2}\omega _{1}-\left. P_{\rho }^{-1}\pi
_{0}\right\vert _{\Omega }=\omega _{2}^{\ast }+\lambda \omega _{1}^{\ast }%
\text{ \ in }\Omega \\
\\
\left. \omega _{1}\right\vert _{\partial \Omega }=\left. \frac{\partial
\omega _{1}}{\partial n}\right\vert _{\partial \Omega }=0\text{ \ on }%
\partial \Omega%
\end{array}%
\right.  \label{s} \\
&  \notag \\
& \left\{
\begin{array}{l}
\lambda \mu -\Delta \mu +\nabla \pi _{0}=\mu ^{\ast }\text{ \ in }\mathcal{O}
\\
\\
\text{div}(\mu )=0\text{ \ in }\mathcal{O} \\
\\
\mu =0\text{ on }S\text{; \ }\mu =[0,0,\lambda \omega _{1}-\omega _{1}^{\ast
}]\text{ \ in }\Omega .%
\end{array}%
\right.  \label{f}
\end{align}%
Here, the pressure term is given, via (\ref{A})-(\ref{domain}), as
\begin{equation}
\pi _{0}=G_{\rho ,1}(\omega _{1})+G_{\rho ,2}(\mu )\in H^{1}(\mathcal{O}),%
\text{ for }\left[ \omega _{1},\omega _{2},\mu \right] \in D(\mathcal{A}%
_{\rho }),  \label{press}
\end{equation}%
where $G_{\rho ,1}$ and $G_{\rho ,2}$ are given by (\ref{G1}) and (\ref{G2}).

\bigskip

The proof of Theorem \ref{back} will ultimately depend on the appropriate
use of four basic relations:

(i) Taking the $D(P_{\rho }^{\frac{1}{2}})$-inner product of both sides of
the structural PDE in (\ref{s}) with $\omega _{1}$, integrating by parts and
subsequently taking the real part of the result, we have
\begin{equation}
\alpha ^{2}(1-\tan ^{2}\vartheta )\left\Vert P_{\rho }^{\frac{1}{2}}\omega
_{1}\right\Vert _{L^{2}(\Omega )}^{2}+\left\Vert \Delta \omega
_{1}\right\Vert _{L^{2}(\Omega )}^{2}=\text{Re}\left( \left. \pi
_{0}\right\vert _{\Omega },\omega _{1}\right) _{\Omega }+\text{Re}\left(
P_{\rho }[\omega _{2}^{\ast }+\lambda \omega _{1}^{\ast }],\omega
_{1}\right) _{\Omega },  \label{p1}
\end{equation}%
(after also using implicitly the Criterion 1 above).

(ii) We take the $\mathbf{L}^{2}(\mathcal{O)}$-inner product of both sides
of the fluid PDE in (\ref{f}) with $\mu $. After integrating by parts and
then taking the respective imaginary and real parts of the resulting
relation, we have,
\begin{align}
\beta \left\Vert \mu \right\Vert _{\mathcal{O}}^{2}& =-\text{Im}\left(
\left. \pi _{0}\right\vert _{\Omega },\lambda \omega _{1}-\omega _{1}^{\ast
}\right) _{\Omega }+\text{Im}\left( \mu ^{\ast },\mu \right) _{\mathcal{O}};
\label{p1.2} \\
&  \notag \\
\alpha \left\Vert \mu \right\Vert _{\mathcal{O}}^{2}+\left\Vert \nabla \mu
\right\Vert _{\mathcal{O}}^{2}& =-\text{Re}\left( \left. \pi _{0}\right\vert
_{\Omega },\lambda \omega _{1}-\omega _{1}^{\ast }\right) _{\Omega }+\text{Re%
}\left( \mu ^{\ast },\mu \right) _{\mathcal{O}}.  \label{1.22}
\end{align}

(iii) Lastly, we take the $\mathbf{H}_{\rho }$-inner product of both sides
of the resolvent equation (\ref{abstract}) with respect to solution
variables $[\omega _{1},\omega _{2},\mu ]$. This gives, upon integrating and
taking the real part of the resulting relation:%
\begin{equation}
\alpha \left\Vert \left[
\begin{array}{c}
\omega _{1} \\
\omega _{2} \\
\mu%
\end{array}%
\right] \right\Vert _{\mathbf{H}_{\rho }}^{2}=-\left\Vert \nabla \mu
\right\Vert _{\mathcal{O}}^{2}+\text{Re}\left( \left[
\begin{array}{c}
\omega _{1}^{\ast } \\
\omega _{2}^{\ast } \\
\mu ^{\ast }%
\end{array}%
\right] ,\left[
\begin{array}{c}
\omega _{1} \\
\omega _{2} \\
\mu%
\end{array}%
\right] \right) _{\mathbf{H}_{\rho }}.  \label{p1.3}
\end{equation}

In view of the right hand side of the relations (\ref{p1}) and (\ref{p1.2})-(%
\ref{1.22}), it is evidently necessary to scrutinize the \textquotedblleft
interface\textquotedblright\ term $\left( \left. \pi _{0}\right\vert
_{\Omega },\omega _{1}\right) _{\Omega }$. Indeed, the estimation of this
term will constitute the bulk of the effort in this work. By way of
attaining a useful estimation, we will need to consider the explicit
representation of the pressure term $\pi _{0}$, as given in (\ref{press}).
Via this expression we have then,%
\begin{equation}
\left( \left. \pi _{0}\right\vert _{\Omega },\omega _{1}\right) _{\Omega
}=\left( \left. G_{\rho ,1}(\omega _{1})\right\vert _{\Omega },\omega
_{1}\right) _{\Omega }+\left( \left. G_{\rho ,2}(\mu )\right\vert _{\Omega
},\omega _{1}\right) _{\Omega }.  \label{p2}
\end{equation}%
We will proceed now to estimate each inner product on the right hand side of
(\ref{p2}).

\bigskip

\subsection{Analysis of the Term $\left( \left. G_{\protect\rho ,2}(\protect%
\mu )\right\vert _{\Omega },\protect\omega _{1}\right) _{\Omega }$}

We recall from (\ref{G2}) that
\begin{equation}
G_{\rho ,2}(\mu )=R_{\rho }(\left. \Delta \mu ^{3}\right\vert _{\Omega })+%
\tilde{R}_{p}(\left. \Delta \mu \cdot \nu \right\vert _{S}).  \label{2.1}
\end{equation}%
With the right hand side of (\ref{2.1}) in mind we define the positive,
self-adjoint operator $B_{\rho }:D(B_{\rho })\subset L^{2}(\mathcal{O}%
)\rightarrow L^{2}(\mathcal{O})$ by%
\begin{equation}
B_{\rho }f=-\Delta f\text{ \ in }\mathcal{O};\,\,D(B_{\rho })=\left\{ f\in
H^{1}(\mathcal{O}):\Delta f\in L^{2}(\mathcal{O})\text{ and }\left\{
\begin{array}{l}
\frac{\partial f}{\partial \nu }+P_{\rho }^{-1}f=0\text{ \ on }\Omega \\
\frac{\partial f}{\partial \nu }=0\text{ \ on }S%
\end{array}%
\right. \text{ \ }\right\} .  \label{B}
\end{equation}

Therewith one can can readily compute the respective adjoints of $R_{\rho
}\in \mathcal{L}(H^{-\frac{1}{2}}(\Omega ),H^{1}(\mathcal{O}))$, $\tilde{R}%
_{\rho }\in \mathcal{L}(H^{-\frac{1}{2}}(S),H^{1}(\mathcal{O}))$, $B_{\rho
}R_{\rho }\in \mathcal{L}(H^{-\frac{1}{2}}(\Omega ),[H^{1}(\mathcal{O}%
)]^{\prime })$ and $B_{\rho }\tilde{R}_{\rho }\in \mathcal{L}(H^{-\frac{1}{2}%
}(S),[H^{1}(\mathcal{O})]^{\prime })$, as,%
\begin{align}
R_{\rho }^{\ast }f& =\left. B_{\rho }^{-1}f\right\vert _{\Omega }\text{ for
all }f\in \lbrack H^{1}(\mathcal{O})]^{\prime };  \label{2.2} \\
\tilde{R}_{\rho }^{\ast }f& =\left. B_{\rho }^{-1}f\right\vert _{S}\text{
for all }f\in \lbrack H^{1}(\mathcal{O})]^{\prime };  \label{2.3} \\
R_{\rho }^{\ast }B_{\rho }f& =\left. f\right\vert _{\Omega }\text{ for all }%
f\in \lbrack H^{1}(\mathcal{O})];  \label{2.4} \\
\tilde{R}_{\rho }^{\ast }B_{\rho }f& =\left. f\right\vert _{S}\text{ for all
}f\in \lbrack H^{1}(\mathcal{O})].  \label{2.5}
\end{align}

\bigskip

Indeed, to show (\ref{2.2}): Given $g\in H^{-\frac{1}{2}}(\Omega )$ and $%
f\in \lbrack H^{1}(\mathcal{O})]^{\prime }$, we have from (\ref{B})%
\begin{align*}
(R_{\rho }g,f)_{\mathcal{O}}& =(R_{\rho }g,(-\Delta )B_{\rho }^{-1}f)_{%
\mathcal{O}} \\
& =(\nabla R_{\rho }g,\nabla B_{\rho }^{-1}f)_{\mathcal{O}}-(R_{\rho }g,%
\frac{\partial }{\partial \nu }B_{\rho }^{-1}f)_{\Omega }+0 \\
& =(-\Delta R_{\rho }g,B_{\rho }^{-1}f)_{\mathcal{O}}+(\frac{\partial }{%
\partial \nu }R_{\rho }g,B_{\rho }^{-1}f)_{\Omega }+(R_{\rho }g,P_{\rho
}^{-1}B_{\rho }^{-1}f)_{\Omega } \\
& =(\frac{\partial }{\partial \nu }R_{\rho }g+P_{\rho }^{-1}R_{\rho
}g,B_{\rho }^{-1}f)_{\Omega } \\
& =(g,B_{\rho }^{-1}f)_{\Omega }.
\end{align*}%
The proofs of relations (\ref{2.3})-(\ref{2.5}) are similar.

\bigskip

With the relations (\ref{2.2})-(\ref{2.5}) in hand, we proceed: From (\ref%
{2.1}) we have%
\begin{align*}
\left( \left. G_{\rho ,2}(\mu )\right\vert _{\Omega },\omega _{1}\right)
_{\Omega }& =\left( R_{\rho }^{\ast }B_{\rho }\left[ R_{\rho }(\left. \Delta
\mu ^{3}\right\vert _{\Omega })+\tilde{R}_{p}(\left. \Delta \mu \cdot \nu
\right\vert _{S})\right] ,\omega _{1}\right) _{\Omega } \\
& =\left( \left[ R_{\rho }(\left. \Delta \mu ^{3}\right\vert _{\Omega })+%
\tilde{R}_{p}(\left. \Delta \mu \cdot \nu \right\vert _{S})\right] ,B_{\rho
}R_{\rho }\omega _{1}\right) _{\mathcal{O}} \\
& =\left( \left. \Delta \mu ^{3}\right\vert _{\Omega },\left. R_{\rho
}\omega _{1}\right\vert _{\Omega }\right) _{\Omega }+\left( \left. \Delta
\mu \cdot \nu \right\vert _{S},\left. R_{\rho }\omega _{1}\right\vert
_{S}\right) _{S} \\
& =\left( \Delta \mu \cdot \nu ,R_{\rho }\omega _{1}\right) _{\partial
\mathcal{O}}
\end{align*}%
(and we are also using here the fact that $\left. \nu \right\vert _{\Omega
}=[0,0,1]$). Invoking now Green's Formula -- and simultaneously using the
fact that fluid term $\Delta \mu $ is divergence free -- yields%
\begin{equation}
\left( \left. G_{\rho ,2}(\mu )\right\vert _{\Omega },\omega _{1}\right)
_{\Omega }=\left( \Delta \mu ,\nabla R_{\rho }\omega _{1}\right) _{\mathcal{O%
}}\text{.}  \label{2.7}
\end{equation}%
Following this relation up with Green's First Identity, we have then%
\begin{equation}
\left( \left. G_{\rho ,2}(\mu )\right\vert _{\Omega },\omega _{1}\right)
_{\Omega }=-\left( \nabla \mu ,\nabla (\nabla R_{\rho }\omega _{1})\right) _{%
\mathcal{O}}+\left\langle \frac{\partial \mu }{\partial \nu },\left. \nabla
R_{\rho }\omega _{1}\right\vert _{\partial \mathcal{O}}\right\rangle
_{\partial \mathcal{O}}\text{.}  \label{2.8}
\end{equation}

\bigskip

\subsubsection{Estimating the term $\left\vert \left( \protect\nabla \protect%
\mu ,\protect\nabla (\protect\nabla R_{\protect\rho }\protect\omega %
_{1})\right) _{\mathcal{O}}\right\vert $}

At this point, we consider the term $R_{\rho }\omega _{1}$ -- where map $%
R_{\rho }$ is again given by (\ref{R1}) -- as the solution of the following
elliptic problem: For $\rho \geq 0$, we see from (\ref{R1}) that $R_{\rho
}\omega _{1}$ solves the elliptic problem%
\begin{eqnarray}
\Delta (R_{\rho }\omega _{1}) &=&0\text{ \ in }\mathcal{O}  \notag \\
\frac{\partial (R_{\rho }\omega _{1})}{\partial \nu } &=&\left[ \omega
_{1}-P_{\rho }^{-1}(\left. R_{\rho }\omega _{1}\right\vert _{\Omega })\right]
_{ext}\text{ \ on }\partial \mathcal{O}\text{,}  \label{bvp}
\end{eqnarray}%
where $L^{2}(\partial \mathcal{O})$-Neumann data is given by%
\begin{equation*}
\left[ \omega _{1}-P_{\rho }^{-1}(\left. R_{\rho }\omega _{1}\right\vert
_{\Omega })\right] _{ext}\equiv \left\{
\begin{array}{l}
\omega _{1}-P_{\rho }^{-1}(\left. R_{\rho }\omega _{1}\right\vert _{\Omega })%
\text{ \ on }\Omega \\
0\text{ \ on\ }S.%
\end{array}%
\right.
\end{equation*}%
Then by the regularity result in \cite{jerison}, valid for Lipschitz
domains, we have the estimate%
\begin{eqnarray}
\left\Vert R_{\rho }\omega _{1}\right\Vert _{H^{\frac{3}{2}}(\mathcal{O})}
&\leq &C\left\Vert \left[ \omega _{1}-P_{\rho }^{-1}(\left. R_{\rho }\omega
_{1}\right\vert _{\Omega })\right] _{ext}\right\Vert _{L^{2}(\partial
\mathcal{O})}  \notag \\
&\leq &C\left( \left\Vert \omega _{1}\right\Vert _{L^{2}(\Omega
)}+\left\Vert P_{\rho }^{-1}(\left. R_{\rho }\omega _{1}\right\vert _{\Omega
})\right\Vert _{L^{2}(\Omega )}\right)  \notag \\
&\leq &C\left\Vert \omega _{1}\right\Vert _{L^{2}(\Omega )},  \label{armed}
\end{eqnarray}%
where in the second to last inequality we have used $P_{\rho }^{-1}\in
\mathcal{L}\left( L^{2}(\Omega ),D(P_{\rho })\right) $, as well as the
boundedness of $R_{\rho }\in \mathcal{L}\left( H^{-\frac{1}{2}}(\Omega ),%
\mathbf{H}^{1}(\mathcal{O})\right) $, which is noted in (\ref{Rs}).

Using the estimate in (\ref{armed}), in tandem with interpolation, we have
now -- using implicitly $H^{\frac{1}{2}}(\mathcal{O})=H_{0}^{\frac{1}{2}}(%
\mathcal{O})$; see e.g., Theorem 3.40 (i) of \cite{mclean} --%
\begin{align*}
\left( \nabla \mu ,\nabla (\nabla R_{\rho }\omega _{1})\right) _{\mathcal{O}%
}& =\left\langle \nabla \mu ,\nabla (\nabla R_{\rho }\omega
_{1})\right\rangle _{H^{\frac{1}{2}}(\mathcal{O})\times H^{-\frac{1}{2}}(%
\mathcal{O})} \\
& \leq C\left\Vert \nabla \mu \right\Vert _{H^{\frac{1}{2}}(\mathcal{O}%
)}\left\Vert \nabla (\nabla R_{\rho }\omega _{1})\right\Vert _{H^{-\frac{1}{2%
}}(\mathcal{O})} \\
& \leq C\left\Vert \mu \right\Vert _{H^{\frac{3}{2}}(\mathcal{O})}\left\Vert
\omega _{1}\right\Vert _{L^{2}(\Omega )} \\
& \leq C\left\Vert \mu \right\Vert _{H^{1}(\mathcal{O})}^{\frac{1}{2}%
}\left\Vert \mu \right\Vert _{H^{2}(\mathcal{O})}^{\frac{1}{2}}\left\Vert
\omega _{1}\right\Vert _{L^{2}(\Omega )} \\
& \leq C\left\Vert \mu \right\Vert _{H^{1}(\mathcal{O})}^{\frac{1}{2}%
}\left\Vert [\omega _{1},\omega _{2},\mu ]\right\Vert _{D(\mathcal{A}_{\rho
})}^{\frac{1}{2}}\left\Vert \omega _{1}\right\Vert _{L^{2}(\Omega )} \\
& =C\frac{\left\vert \alpha \right\vert ^{\frac{1}{4}}}{\left\vert \alpha
\right\vert ^{\frac{1}{4}}}\left\Vert \mu \right\Vert _{H^{1}(\mathcal{O})}^{%
\frac{1}{2}}\left\Vert [\omega _{1},\omega _{2},\mu ]\right\Vert _{D(%
\mathcal{A}_{\rho })}^{\frac{1}{2}}\left\Vert \omega _{1}\right\Vert
_{L^{2}(\Omega )}.
\end{align*}%
Estimating further the right hand side of this inequality, via Young's
Inequality, now yields%
\begin{align}
\left\vert \left( \nabla \mu ,\nabla (\nabla R_{\rho }\omega _{1})\right) _{%
\mathcal{O}}\right\vert & \leq \frac{\epsilon }{\left\vert \alpha
\right\vert }\left\Vert \nabla \mu \right\Vert _{\mathcal{O}%
}^{2}+C_{\epsilon }\left\vert \alpha \right\vert ^{\frac{1}{3}}\left\Vert
[\omega _{1},\omega _{2},\mu ]\right\Vert _{D(\mathcal{A}_{\rho })}^{\frac{2%
}{3}}\left\Vert \omega _{1}\right\Vert _{L^{2}(\Omega )}^{\frac{4}{3}}
\notag \\
& =\frac{\epsilon }{\left\vert \alpha \right\vert }\left\Vert \nabla \mu
\right\Vert _{\mathcal{O}}^{2}+C_{\epsilon }\left\vert \alpha \right\vert
\left( 1+\tan ^{2}\vartheta \right) ^{\frac{1}{3}}\left\Vert \left[
\begin{array}{c}
\omega _{1} \\
\omega _{2} \\
\mu%
\end{array}%
\right] -\frac{1}{\left\vert \alpha \right\vert e^{\pm i\vartheta }\sqrt{%
1+\tan ^{2}\vartheta }}\left[
\begin{array}{c}
\omega _{1}^{\ast } \\
\omega _{2}^{\ast } \\
\mu ^{\ast }%
\end{array}%
\right] \right\Vert _{\mathbf{H}_{\rho }}^{\frac{2}{3}}  \notag \\
& \text{ \ \ \ \ }\times \left\Vert \omega _{1}\right\Vert _{L^{2}(\Omega
)}^{\frac{4}{3}},  \label{2.9a}
\end{align}%
where in the last step, we have used the resolvent relation (\ref{abstract}%
), and the assumption in Criterion 1 that frequency domain parameter $%
\lambda $ lies along one of the two rays $e^{\pm i\vartheta }$, for fixed $%
\vartheta \in (3\pi /4,\pi )$. Estimating once more, we have then%
\begin{eqnarray*}
\left\vert \left( \nabla \mu ,\nabla (\nabla R_{\rho }\omega _{1})\right) _{%
\mathcal{O}}\right\vert &\leq &\frac{\epsilon }{\left\vert \alpha
\right\vert }\left\Vert \nabla \mu \right\Vert _{\mathcal{O}}^{2}+\frac{%
\epsilon }{2}\left\Vert \left[
\begin{array}{c}
\omega _{1} \\
\omega _{2} \\
\mu%
\end{array}%
\right] -\frac{1}{\left\vert \alpha \right\vert e^{\pm i\vartheta }\sqrt{%
1+\tan ^{2}\vartheta }}\left[
\begin{array}{c}
\omega _{1}^{\ast } \\
\omega _{2}^{\ast } \\
\mu ^{\ast }%
\end{array}%
\right] \right\Vert _{\mathbf{H}_{\rho }}^{2} \\
&&+C_{\epsilon ,\vartheta ,\delta }\left\vert \alpha \right\vert ^{\frac{3}{2%
}}\left\Vert \omega _{1}\right\Vert _{\Omega }^{2},
\end{eqnarray*}%
or for $\left\vert \alpha \right\vert >1,$%
\begin{equation}
\left\vert \left( \nabla \mu ,\nabla (\nabla R_{\rho }\omega _{1})\right) _{%
\mathcal{O}}\right\vert \leq \frac{\epsilon }{\left\vert \alpha \right\vert }%
\left\Vert \nabla \mu \right\Vert _{\mathcal{O}}^{2}+\epsilon \left\Vert %
\left[
\begin{array}{c}
\omega _{1} \\
\omega _{2} \\
\mu%
\end{array}%
\right] \right\Vert _{\mathbf{H}_{\rho }}^{2}+C_{\epsilon ,\vartheta ,\delta
}\left\vert \alpha \right\vert ^{\frac{3}{2}}\left\Vert \omega
_{1}\right\Vert _{D(P_{\rho }^{\frac{1}{2}})}^{2}+C_{\epsilon ,\delta
}\left\Vert \left[
\begin{array}{c}
\omega _{1}^{\ast } \\
\omega _{2}^{\ast } \\
\mu ^{\ast }%
\end{array}%
\right] \right\Vert _{\mathbf{H}_{\rho }}^{2}.  \label{2.9}
\end{equation}

\bigskip

\subsubsection{Estimating the term $\left\vert \left\langle \frac{\partial
\protect\mu }{\partial \protect\nu },\left. \protect\nabla R_{\protect\rho }%
\protect\omega _{1}\right\vert _{\partial \mathcal{O}}\right\rangle
_{\partial \mathcal{O}}\right\vert $ for $\protect\rho =0$.}

The second term on the right hand side of (\ref{2.8}) is an even more
delicate matter; the analysis here necessarily becomes a dichotomy with
respect to $\rho =0$ and $\rho >0$. In either case, we will need the
following boundary trace inequalities (see e.g., Theorem 1.6.6 of \cite%
{brenner}, p. 37): Let $D$ be a bounded domain in $\mathbb{R}^{n}$, $n\geq 2$%
, with Lipschitz boundary $\partial D$. Then there is a positive constant $%
C^{\ast }$such that
\begin{align}
\left\Vert f\right\Vert _{\partial D}& \leq C^{\ast }\left\Vert f\right\Vert
_{L^{2}(D)}^{\frac{1}{2}}\left\Vert f\right\Vert _{H^{1}(D\mathcal{)}}^{%
\frac{1}{2}}\text{ for every }f\in H^{1}(D).  \label{trace} \\
\left\Vert \frac{\partial f}{\partial \nu }\right\Vert _{\partial D}& \leq
C^{\ast }\left\Vert f\right\Vert _{H^{1}(D\mathcal{)}}^{\frac{1}{2}%
}\left\Vert f\right\Vert _{H^{2}(D\mathcal{)}}^{\frac{1}{2}}\text{ for every
}f\in H^{2}(D).  \label{trace2}
\end{align}%
Note that the second inequality follows from (\ref{trace}), after using the
fact that normal vector $\nu \in \mathbf{L}^{\infty }(\partial D)$, since $%
\partial D$ is Lipschitz; see \cite{necas} (and so constant $C^{\ast }$
depends upon $\left\Vert \nu \right\Vert _{\mathbf{L}^{\infty }(\partial D)}$%
.)

\bigskip

To start: We will have need here of the following positive definite,
self-adjoint operator $\mathbf{\mathring{A}}$ $:D(\mathring{A})\subset
L^{2}(\Omega )\rightarrow L^{2}(\Omega )$, defined by%
\begin{equation}
\mathbf{\mathring{A}}\varpi =\Delta ^{2}\varpi \text{, \ }D(\mathbf{%
\mathring{A}})=H^{4}(\Omega )\cap H_{0}^{2}(\Omega ).  \label{ang}
\end{equation}%
(Note that in the case that $\Omega $ is polygonal -- i.e., geometric
condition \textsf{(G.2)} is in force -- the angle condition assumed in
\textsf{(G.2)} assures the smoothness of $D(\mathbf{\mathring{A}})$ as
given; see Theorem 2 of \cite{blum}). As such, this operator obeys the
following \textquotedblleft analyticity\textquotedblright\ estimate for all $%
s>0:$%
\begin{equation}
\left\Vert \mathbf{\mathring{A}}^{\eta }\mathcal{R}(-s;\mathbf{\mathring{A}}%
)\right\Vert _{\mathcal{L}(L^{2}(\Omega ))}\leq \frac{C}{(1+s)^{1-\eta }}%
\text{, for all }\eta \in \lbrack 0,1]  \label{est}
\end{equation}%
(see e.g., the expression (5.15) in \cite{krein}, p. 115). With this
operator in hand, then in the present case $\rho =0$ the structural equation
in (\ref{s}) can be written as%
\begin{equation*}
\left[ \alpha ^{2}(1-\tan ^{2}\vartheta )+\mathbf{\mathring{A}}\right]
\omega _{1}=-2i\alpha \beta \omega _{1}+\left. \pi _{0}\right\vert _{\Omega
}+\omega _{2}^{\ast }+\lambda \omega _{1}^{\ast }
\end{equation*}%
(after also using Criterion 1). Applying thereto the operator $\mathbf{%
\mathring{A}}^{\eta }\mathcal{R}(-\alpha ^{2}(1-\tan ^{2}\vartheta );\mathbf{%
\mathring{A}})$ gives then
\begin{equation*}
\mathbf{\mathring{A}}^{\eta }\omega _{1}=\mathbf{\mathring{A}}^{\eta }%
\mathcal{R}(-\alpha ^{2}(1-\tan ^{2}\vartheta );\mathbf{\mathring{A}})\left[
2i\alpha \beta \omega _{1}-\left. \pi _{0}\right\vert _{\Omega }-\omega
_{2}^{\ast }-\lambda \omega _{1}^{\ast }\right] .
\end{equation*}%
Subsequently applying the estimate (\ref{est}), we then have for $0\leq \eta
\leq 1$ and $\left\vert \alpha \right\vert >0$ sufficently large,
\begin{align}
\left\Vert \mathbf{\mathring{A}}^{\eta }\omega _{1}\right\Vert _{\Omega }&
\leq \frac{C}{(1+\alpha ^{2}(1-\tan ^{2}\vartheta ))^{1-\eta }}\left[
\left\vert \alpha \beta \right\vert \left\Vert \omega _{1}\right\Vert
_{L^{2}(\Omega )}+\left\Vert \left. \pi _{0}\right\vert _{\Omega
}\right\Vert _{L^{2}(\Omega )}+\left\Vert \omega _{2}^{\ast }+\lambda \omega
_{1}^{\ast }\right\Vert _{L^{2}(\Omega )}\right]  \notag \\
& \leq C_{\vartheta }\left\vert \alpha \right\vert ^{2\eta }\left\Vert
\omega _{1}\right\Vert _{\Omega }+\frac{C_{\vartheta }}{\left\vert \alpha
\right\vert ^{1-2\eta }}\left( \left\Vert [\omega _{1},\omega _{2},\mu
]\right\Vert _{\mathbf{H}_{0}}+\left\Vert [\omega _{1}^{\ast },\omega
_{2}^{\ast },\mu ^{\ast }]\right\Vert _{\mathbf{H}_{0}}\right) .  \label{eta}
\end{align}%
In obtaining this estimate, we have used $\left\vert \beta \right\vert
=\left\vert \alpha \right\vert \left\vert \tan \vartheta \right\vert $, the
expression (\ref{press}), and the resolvent equation (\ref{abstract}).

\bigskip

With estimate (\ref{eta}) in hand, we now estimate the second term on the
right hand side of (\ref{2.8}): Reinvoking the estimate (\ref{armed}) for
the solution of \ref{bvp}) (with therein $\rho =0$), in combination with the
trace inequality (\ref{trace2}), we have for $\left\vert \alpha \right\vert
>0$ sufficiently large,
\begin{align*}
\left\vert \left\langle \frac{\partial \mu }{\partial \nu },\left. \nabla
R_{0}\omega _{1}\right\vert _{\partial \mathcal{O}}\right\rangle _{\partial
\mathcal{O}}\right\vert & \leq \left\Vert \frac{\partial \mu }{\partial \nu }%
\right\Vert _{\mathbf{L}^{2}(\partial \mathcal{O})}\left\Vert \left. \nabla
R_{0}\omega _{1}\right\vert _{\partial \mathcal{O}}\right\Vert
_{L^{2}(\partial \mathcal{O})} \\
& \leq C\left\Vert \nabla \mu \right\Vert _{\mathbf{L}^{2}(\mathcal{O})}^{%
\frac{1}{2}}\left\Vert \mu \right\Vert _{\mathbf{H}^{2}(\mathcal{O})}^{\frac{%
1}{2}}\left\Vert \omega _{1}\right\Vert _{L^{2}(\Omega \mathcal{)}} \\
& \leq C\left\Vert \nabla \mu \right\Vert _{\mathbf{L}^{2}(\mathcal{O})}^{%
\frac{1}{2}}\left\Vert [\omega _{1},\omega _{2},\mu ]\right\Vert _{D(%
\mathcal{A}_{0})}^{\frac{1}{2}}\left\Vert \omega _{1}\right\Vert
_{L^{2}(\Omega )} \\
& =C_{\vartheta }\sqrt{\left\vert \alpha \right\vert }\left\Vert \nabla \mu
\right\Vert _{\mathbf{L}^{2}(\mathcal{O})}^{\frac{1}{2}}\left\Vert \left[
\begin{array}{c}
\omega _{1} \\
\omega _{2} \\
\mu%
\end{array}%
\right] -\frac{1}{\left\vert \alpha \right\vert e^{\pm i\vartheta }\sqrt{%
1+\tan ^{2}\vartheta }}\left[
\begin{array}{c}
\omega _{1}^{\ast } \\
\omega _{2}^{\ast } \\
\mu ^{\ast }%
\end{array}%
\right] \right\Vert _{\mathbf{H}_{0}}^{\frac{1}{2}}\left\Vert \omega
_{1}\right\Vert _{L^{2}(\Omega )},
\end{align*}%
where again we have implicitly used the resolvent relation (\ref{abstract}).
Using now the characterization
\begin{equation*}
H^{\frac{1}{2}-\delta }(\Omega )\approx D(\mathbf{\mathring{A}}^{\frac{1}{8}-%
\frac{\delta }{4}}),
\end{equation*}%
which can be inferred from the definition of the domain in (\ref{ang}) and
\cite{grisvard}, \ we have upon applying (\ref{eta}) to the right hand side
of (\ref{eta2}),%
\begin{align*}
\left\vert \left\langle \frac{\partial \mu }{\partial \nu },\left. \nabla
R_{0}\omega _{1}\right\vert _{\partial \mathcal{O}}\right\rangle _{\partial
\mathcal{O}}\right\vert & \leq C_{\vartheta }\sqrt{\left\vert \alpha
\right\vert }\left\Vert \nabla \mu \right\Vert _{\mathcal{O}}^{\frac{1}{2}%
}\left\Vert \left[
\begin{array}{c}
\omega _{1} \\
\omega _{2} \\
\mu%
\end{array}%
\right] -\frac{1}{\left\vert \alpha \right\vert e^{\pm i\vartheta }\sqrt{%
1+\tan ^{2}\vartheta }}\left[
\begin{array}{c}
\omega _{1}^{\ast } \\
\omega _{2}^{\ast } \\
\mu ^{\ast }%
\end{array}%
\right] \right\Vert _{\mathbf{H}_{0}}^{\frac{1}{2}}\times \\
& \text{ \ \ \ \ \ \ \ \ }\left[ \left\vert \alpha \right\vert ^{\frac{1}{4}-%
\frac{\delta }{2}}\left\Vert \omega _{1}\right\Vert _{\Omega }+\frac{1}{%
\left\vert \alpha \right\vert ^{\frac{3}{4}+\frac{\delta }{2}}}\left(
\left\Vert [\omega _{1},\omega _{2},\mu ]\right\Vert _{\mathbf{H}%
_{0}}+\left\Vert [\omega _{1}^{\ast },\omega _{2}^{\ast },\mu ^{\ast
}]\right\Vert _{\mathbf{H}_{0}}\right) \right] \\
& \\
& =\frac{C_{\vartheta }}{\left\vert \alpha \right\vert ^{\frac{1}{4}}}%
\left\Vert \nabla \mu \right\Vert _{\mathcal{O}}^{\frac{1}{2}}\left\vert
\alpha \right\vert ^{1-\frac{\delta }{2}}\left\Vert \left[
\begin{array}{c}
\omega _{1} \\
\omega _{2} \\
\mu%
\end{array}%
\right] -\frac{1}{\left\vert \alpha \right\vert e^{\pm i\vartheta }\sqrt{%
1+\tan ^{2}\vartheta }}\left[
\begin{array}{c}
\omega _{1}^{\ast } \\
\omega _{2}^{\ast } \\
\mu ^{\ast }%
\end{array}%
\right] \right\Vert _{\mathbf{H}_{0}}^{\frac{1}{2}}\times \\
& \text{ \ \ \ \ \ \ }\left[ \left\Vert \omega _{1}\right\Vert _{\Omega }+%
\frac{1}{\left\vert \alpha \right\vert }\left( \left\Vert [\omega
_{1},\omega _{2},\mu ]\right\Vert _{\mathbf{H}_{0}}+\left\Vert [\omega
_{1}^{\ast },\omega _{2}^{\ast },\mu ^{\ast }]\right\Vert _{\mathbf{H}%
_{0}}\right) \right] .
\end{align*}%
This gives now, via Young's Inequality, for $\left\vert \alpha \right\vert
>0 $ sufficiently large,%
\begin{align}
\left\vert \left\langle \frac{\partial \mu }{\partial \nu },\left. \nabla
R_{0}\omega _{1}\right\vert _{\partial \mathcal{O}}\right\rangle _{\partial
\mathcal{O}}\right\vert & \leq \frac{\epsilon }{\left\vert \alpha
\right\vert }\left\Vert \nabla \mu \right\Vert _{\mathcal{O}%
}^{2}+C_{\epsilon ,\vartheta }\left\Vert \left[
\begin{array}{c}
\omega _{1} \\
\omega _{2} \\
\mu%
\end{array}%
\right] -\frac{1}{\left\vert \alpha \right\vert e^{\pm i\vartheta }\sqrt{%
1+\tan ^{2}\vartheta }}\left[
\begin{array}{c}
\omega _{1}^{\ast } \\
\omega _{2}^{\ast } \\
\mu ^{\ast }%
\end{array}%
\right] \right\Vert _{\mathbf{H}_{0}}^{\frac{2}{3}}\times  \notag \\
& \,\,\,\,\,\text{\ }\left\vert \alpha \right\vert ^{\frac{4}{3}-\frac{%
2\delta }{3}}\left[ \left\Vert \omega _{1}\right\Vert _{\Omega }+\frac{1}{%
\left\vert \alpha \right\vert }\left( \left\Vert [\omega _{1},\omega
_{2},\mu ]\right\Vert _{\mathbf{H}_{0}}+\left\Vert [\omega _{1}^{\ast
},\omega _{2}^{\ast },\mu ^{\ast }]\right\Vert _{\mathbf{H}_{0}}\right) %
\right] ^{\frac{4}{3}}  \notag \\
&  \notag \\
& \leq \frac{\epsilon }{\left\vert \alpha \right\vert }\left\Vert \nabla \mu
\right\Vert _{\mathcal{O}}^{2}+\epsilon \left\Vert \lbrack \omega
_{1},\omega _{2},\mu ]\right\Vert _{\mathbf{H}_{0}}^{2}  \notag \\
& \,\,\,\,\,+C_{\epsilon ,\vartheta }\left\vert \alpha \right\vert
^{2-\delta }\left\Vert \omega _{1}\right\Vert _{\Omega }^{2}+\frac{%
C_{\epsilon ,\vartheta }}{\left\vert \alpha \right\vert ^{\delta }}\left(
\left\Vert [\omega _{1},\omega _{2},\mu ]\right\Vert _{\mathbf{H}%
_{0}}^{2}+\left\Vert [\omega _{1}^{\ast },\omega _{2}^{\ast },\mu ^{\ast
}]\right\Vert _{\mathbf{H}_{0}}^{2}\right) .  \label{eta_f}
\end{align}

\bigskip

\subsubsection{Estimating the term $\left\vert \left\langle \frac{\partial
\protect\mu }{\partial \protect\nu },\left. \protect\nabla R_{\protect\rho }%
\protect\omega _{1}\right\vert _{\partial \mathcal{O}}\right\rangle
_{\partial \mathcal{O}}\right\vert $ for $\protect\rho >0$.}

Using again the estimates in (\ref{armed}) and (\ref{trace2}), along with
the Sobolev Trace Theorem, we have for $0<\delta <\frac{1}{2}$,%
\begin{align*}
\left\vert \left\langle \frac{\partial \mu }{\partial \nu },\left. \nabla
R_{\rho }\omega _{1}\right\vert _{\partial \mathcal{O}}\right\rangle
_{\partial \mathcal{O}}\right\vert & \leq \left\Vert \frac{\partial \mu }{%
\partial \nu }\right\Vert _{\partial \mathcal{O}}\left\Vert \left. \nabla
R_{\rho }\omega _{1}\right\vert _{\partial \mathcal{O}}\right\Vert
_{\partial \mathcal{O}} \\
& \leq C\left\Vert \nabla \mu \right\Vert _{\mathcal{O}}^{\frac{1}{2}%
}\left\Vert \mu \right\Vert _{\mathbf{H}^{2}(\mathcal{O})}^{\frac{1}{2}%
}\left\Vert R_{\rho }\omega _{1}\right\Vert _{H^{2-\delta }(\mathcal{O})} \\
& \leq C\left\Vert \nabla \mu \right\Vert _{\mathcal{O}}^{\frac{1}{2}%
}\left\Vert \mu \right\Vert _{\mathbf{H}^{2}(\mathcal{O})}^{\frac{1}{2}%
}\left\Vert \omega _{1}\right\Vert _{H^{\frac{1}{2}-\delta }(\Omega )} \\
& \leq C\left\Vert \nabla \mu \right\Vert _{\mathcal{O}}^{\frac{1}{2}%
}\left\Vert \mu \right\Vert _{\mathbf{H}^{2}(\mathcal{O})}^{\frac{1}{2}%
}\left\Vert \omega _{1}\right\Vert _{\Omega }^{\delta }\left\Vert \omega
_{1}\right\Vert _{H^{1}(\Omega )}^{1-\delta }.
\end{align*}%
Combining this with the fluid boundary condition in (\ref{f}) and the
resolvent relation (\ref{abstract}), we have then for $\left\vert \alpha
\right\vert >0$ sufficiently large,%
\begin{eqnarray*}
&&\left\vert \left\langle \frac{\partial \mu }{\partial \nu },\left. \nabla
R_{\rho }\omega _{1}\right\vert _{\partial \mathcal{O}}\right\rangle
_{\partial \mathcal{O}}\right\vert \\
&\leq &C\left\Vert \nabla \mu \right\Vert _{\mathcal{O}}^{\frac{1}{2}%
}\left\Vert \mu \right\Vert _{\mathbf{H}^{2}(\mathcal{O})}^{\frac{1}{2}%
}\left\Vert \frac{1}{\lambda }\left( \left. \mu ^{3}\right\vert _{\Omega
}+\omega _{1}^{\ast }\right) \right\Vert _{\Omega }^{\delta }\left\Vert
\omega _{1}\right\Vert _{H^{1}(\Omega )}^{1-\delta } \\
&\leq &C\left\Vert \nabla \mu \right\Vert _{\mathcal{O}}^{\frac{1}{2}%
}\left\Vert [\omega _{1},\omega _{2},\mu ]\right\Vert _{D(\mathcal{A}_{\rho
})}^{\frac{1}{2}}\left\Vert \frac{1}{\lambda }\left( \left. \mu
^{3}\right\vert _{\Omega }+\omega _{1}^{\ast }\right) \right\Vert _{\Omega
}^{\delta }\left\Vert \omega _{1}\right\Vert _{H^{1}(\Omega )}^{1-\delta } \\
&\leq &C_{}\left\Vert \nabla \mu \right\Vert _{\mathcal{O}}^{\frac{1}{2}%
}\left\Vert [\omega _{1},\omega _{2},\mu ]-\frac{[\omega _{1}^{\ast },\omega
_{2}^{\ast },\mu ^{\ast }]}{\left\vert \alpha \right\vert e^{\pm i\vartheta }%
\sqrt{1+\tan ^{2}\vartheta }}\right\Vert _{\mathbf{H}_{\rho }}^{\frac{1}{2}%
}\left\Vert \left. \mu ^{3}\right\vert _{\Omega }+\omega _{1}^{\ast
}\right\Vert _{\Omega }^{\delta }\left\Vert \omega _{1}\right\Vert
_{H^{1}(\Omega )}^{1-\delta } \\
&=&C_{ }\frac{\left\vert \alpha \right\vert ^{\frac{1}{4}}}{\left\vert
\alpha \right\vert ^{\frac{1}{4}}}\left\Vert \nabla \mu \right\Vert _{%
\mathcal{O}}^{\frac{1}{2}}\left\Vert [\omega _{1},\omega _{2},\mu ]-\frac{%
[\omega _{1}^{\ast },\omega _{2}^{\ast },\mu ^{\ast }]}{\left\vert \alpha
\right\vert e^{\pm i\vartheta }\sqrt{1+\tan ^{2}\vartheta }}\right\Vert _{%
\mathbf{H}_{\rho }}^{\frac{1}{2}}\left\Vert \left. \mu ^{3}\right\vert
_{\Omega }+\omega _{1}^{\ast }\right\Vert _{\Omega }^{\delta } \left\Vert
\omega _{1}\right\Vert _{H^{1}(\Omega )}^{1-\delta } \\
&\leq &\frac{\epsilon }{\left\vert \alpha \right\vert }\left\Vert \nabla \mu
\right\Vert _{\mathcal{O}}^{2}+ C_{\epsilon} \left\vert \alpha \right\vert ^{%
\frac{1}{3}}\left\Vert [\omega _{1},\omega _{2},\mu ]-\frac{[\omega
_{1}^{\ast },\omega _{2}^{\ast },\mu ^{\ast }]}{\left\vert \alpha
\right\vert e^{\pm i\vartheta }\sqrt{1+\tan ^{2}\vartheta }}\right\Vert _{%
\mathbf{H}_{\rho }}^{\frac{2}{3}}\left\Vert \left. \mu ^{3}\right\vert
_{\Omega }+\omega _{1}^{\ast }\right\Vert _{\Omega }^{\frac{4\delta}{3}%
}\left\Vert \omega _{1}\right\Vert _{H^{1}(\Omega )}^{\frac{4(1-\delta) }{3}}
\\
&\leq &\frac{\epsilon }{\left\vert \alpha \right\vert }\left\Vert \nabla \mu
\right\Vert _{\mathcal{O}}^{2}+\epsilon \left\Vert \lbrack \omega
_{1},\omega _{2},\mu ]\right\Vert _{\mathbf{H}_{\rho }}^{2}+C_{\epsilon
}\left\vert \alpha \right\vert ^{\frac{1}{2}}\left( \frac{\sqrt{\left\vert
\alpha \right\vert }}{\sqrt{\left\vert \alpha \right\vert }}\right) ^{\frac{%
\delta }{2}}\left\Vert \left. \mu ^{3}\right\vert _{\Omega }+\omega
_{1}^{\ast }\right\Vert _{\Omega }^{2\delta }\left\Vert \omega
_{1}\right\Vert _{H^{1}(\Omega )}^{2(1-\delta) } \\
&&\text{ \ \ }+C_{\epsilon }\left\Vert [\omega _{1}^{\ast },\omega
_{2}^{\ast },\mu ^{\ast }]\right\Vert _{\mathbf{H}_{\rho }}^{2} \\
&\leq &\frac{\epsilon }{(C^{\ast })^{2}\left\vert \alpha \right\vert ^{\frac{%
1}{2}}}\left\Vert \left. \mu ^{3}\right\vert _{\Omega }+\omega _{1}^{\ast
}\right\Vert _{\Omega }^{2}+C_{\epsilon }\left\vert \alpha \right\vert ^{%
\frac{1+\delta }{2(1-\delta) }}\left\Vert \omega _{1}\right\Vert
_{H^{1}(\Omega )}^{2}+\frac{\epsilon }{\left\vert \alpha \right\vert }%
\left\Vert \nabla \mu \right\Vert _{\mathcal{O}}^{2}+\epsilon \left\Vert
\lbrack \omega _{1},\omega _{2},\mu ]\right\Vert _{\mathbf{H}_{\rho }}^{2} \\
&&+C_{\epsilon }\left\Vert [\omega _{1}^{\ast },\omega _{2}^{\ast },\mu
^{\ast }]\right\Vert _{\mathbf{H}_{\rho }}^{2},
\end{eqnarray*}%
where $C^{\ast }$ is the positive constant from the interpolation inequality
(\ref{trace}). We have now, for sufficiently large $\left\vert \alpha
\right\vert >1$%
\begin{align}
& \left\vert \left\langle \frac{\partial \mu }{\partial \nu },\left. \nabla
R_{\rho }\omega _{1}\right\vert _{\partial \mathcal{O}}\right\rangle
_{\partial \mathcal{O}}\right\vert  \notag \\
& \leq \frac{\epsilon }{(C^{\ast })^{2}\left\vert \alpha \right\vert ^{\frac{%
1}{2}}}\left\Vert \left. \mu \right\vert _{\partial \mathcal{O}}\right\Vert
_{\partial \mathcal{O}}^{2}+C_{\epsilon }\left\vert \alpha \right\vert ^{%
\frac{1+\delta }{2(1-\delta )}}\left\Vert \omega _{1}\right\Vert
_{H^{1}(\Omega )}^{2}+\frac{\epsilon }{\left\vert \alpha \right\vert }%
\left\Vert \nabla \mu \right\Vert _{\mathcal{O}}^{2}+\epsilon \left\Vert
\lbrack \omega _{1},\omega _{2},\mu ]\right\Vert _{\mathbf{H}_{\rho
}}^{2}+C_{\epsilon }\left\Vert [\omega _{1}^{\ast },\omega _{2}^{\ast },\mu
^{\ast }]\right\Vert _{\mathbf{H}_{\rho }}^{2}  \notag \\
& \leq \frac{\epsilon }{\left\vert \alpha \right\vert ^{\frac{1}{2}}}%
\left\Vert \left. \mu \right\vert \right\Vert _{\mathcal{O}}\left\Vert
\nabla \mu \right\Vert _{\mathcal{O}}+C_{\epsilon }\left\vert \alpha
\right\vert ^{\frac{1+\delta }{2(1-\delta )}}\left\Vert \omega
_{1}\right\Vert _{H^{1}(\Omega )}^{2}+\frac{\epsilon }{\left\vert \alpha
\right\vert }\left\Vert \nabla \mu \right\Vert _{\mathcal{O}}^{2}+\epsilon
\left\Vert \lbrack \omega _{1},\omega _{2},\mu ]\right\Vert _{\mathbf{H}%
_{\rho }}^{2}+C_{\epsilon }\left\Vert [\omega _{1}^{\ast },\omega _{2}^{\ast
},\mu ^{\ast }]\right\Vert _{\mathbf{H}_{\rho }}^{2}  \notag \\
& \leq C_{\epsilon }\left\vert \alpha \right\vert ^{\frac{1+\delta }{%
2(1-\delta )}}\left\Vert \omega _{1}\right\Vert _{H^{1}(\Omega )}^{2}+\frac{%
2\epsilon }{\left\vert \alpha \right\vert }\left\Vert \nabla \mu \right\Vert
_{\mathcal{O}}^{2}+2\epsilon \left\Vert \lbrack \omega _{1},\omega _{2},\mu
]\right\Vert _{\mathbf{H}_{\rho }}^{2}+C_{\epsilon }\left\Vert [\omega
_{1}^{\ast },\omega _{2}^{\ast },\mu ^{\ast }]\right\Vert _{\mathbf{H}_{\rho
}}^{2}.  \label{rho}
\end{align}

\bigskip

Upon a rescaling of small parameter $\epsilon >0$, we have then the estimate
for $\rho >0$,%
\begin{align}
& \left\vert \left\langle \frac{\partial \mu }{\partial \nu },\left. \nabla
R_{\rho }\omega _{1}\right\vert _{\partial \mathcal{O}}\right\rangle
_{\partial \mathcal{O}}\right\vert  \notag \\
& \leq C_{\epsilon }\left\vert \alpha \right\vert ^{\frac{1+\delta }{%
2(1-\delta )}}\left\Vert \omega _{1}\right\Vert _{H^{1}(\Omega )}^{2}+\frac{%
\epsilon }{\left\vert \alpha \right\vert }\left\Vert \nabla \mu \right\Vert
_{\mathcal{O}}^{2}+\epsilon \left\Vert \lbrack \omega _{1},\omega _{2},\mu
]\right\Vert _{\mathbf{H}_{\rho }}^{2}+C_{\epsilon }\left\Vert [\omega
_{1}^{\ast },\omega _{2}^{\ast },\mu ^{\ast }]\right\Vert _{\mathbf{H}_{\rho
}}^{2}.  \label{rho_f}
\end{align}

\bigskip

\bigskip

Combining now (\ref{2.8}), (\ref{2.9}), (\ref{eta_f}), and (\ref{rho_f}),
and taking $\left\vert \alpha \right\vert >0$ sufficently large, we have
finally for all $\rho \geq 0$, and fixed $0<\delta <\frac{1}{2}$,%
\begin{align*}
\left\vert \left( \left. G_{\rho ,2}(\mu )\right\vert _{\Omega },\omega
_{1}\right) _{\Omega }\right\vert & \leq C_{\epsilon, \vartheta }\left\vert
\alpha \right\vert ^{2-\delta }\left\Vert \omega _{1}\right\Vert _{D(P_{\rho
}^{\frac{1}{2}})}^{2}+\frac{\epsilon }{\left\vert \alpha \right\vert }%
\left\Vert \nabla \mu \right\Vert _{\mathcal{O}}^{2} \\
& +\left( \epsilon +\frac{C_{\epsilon , \vartheta}}{\left\vert \alpha
\right\vert ^{\delta }}\right) \left\Vert [\omega _{1},\omega _{2},\mu
]\right\Vert _{\mathbf{H}_{\rho }}^{2}+C_{\epsilon }\left\Vert [\omega
_{1}^{\ast },\omega _{2}^{\ast },\mu ^{\ast }]\right\Vert _{\mathbf{H}_{\rho
}}^{2}.
\end{align*}%
Taking finally $\delta \equiv \frac{1}{2}-\epsilon $, we have then for $%
\left\vert \alpha \right\vert >0$ sufficently large, and $\rho \geq 0$,%
\begin{align}
\left\vert \left( \left. G_{\rho ,2}(\mu )\right\vert _{\Omega },\omega
_{1}\right) _{\Omega }\right\vert & \leq C_{\epsilon, \vartheta }\left\vert
\alpha \right\vert ^{\frac{3}{2}+\epsilon }\left\Vert \omega _{1}\right\Vert
_{D(P_{\rho }^{\frac{1}{2}})}^{2}+\frac{\epsilon }{\left\vert \alpha
\right\vert }\left\Vert \nabla \mu \right\Vert _{\mathcal{O}}^{2}  \notag \\
& +\left( \epsilon +\frac{C_{\epsilon, \vartheta }}{\left\vert \alpha
\right\vert ^{\frac{1}{2}-\epsilon }}\right) \left\Vert [\omega _{1},\omega
_{2},\mu ]\right\Vert _{\mathbf{H}_{\rho }}^{2}+C_{\epsilon, \vartheta
}\left\Vert [\omega _{1}^{\ast },\omega _{2}^{\ast },\mu ^{\ast
}]\right\Vert _{\mathbf{H}\rho }^{2}.  \label{G22}
\end{align}

\bigskip

\subsection{Analysis of the Term $\left\vert \left( \left. G_{\protect\rho %
,1}(\protect\omega _{1})\right\vert _{\Omega },\protect\omega _{1}\right)
_{\Omega }\right\vert $}

\bigskip Recall that the image
\begin{equation}
G_{\rho ,1}(\omega _{1})=R_{\rho }(P_{\rho }^{-1}\Delta ^{2}\omega _{1})
\label{G11}
\end{equation}%
(see (\ref{G1}) and (\ref{p2})).

As before, this work will entail a dichotomy between $\rho =0$ and $\rho >0.$

\subsubsection{Analysis of the term $\left\vert \left( \left. G_{\protect%
\rho ,1}(\protect\omega _{1})\right\vert _{\Omega },\protect\omega %
_{1}\right) _{\Omega }\right\vert $ for $\protect\rho =0$}

In this case, we have from (\ref{G1}) and the expressions in (\ref{2.2}) and
(\ref{2.4}),
\begin{align}
(\left. G_{0,1}(\omega _{1})\right\vert _{\Omega },\omega _{1})_{\Omega }&
=(\left. R_{0}\Delta ^{2}\omega _{1}\right\vert _{\Omega },\omega
_{1})_{\Omega }  \notag \\
& =\left( \Delta ^{2}\omega _{1},\left[ R_{0}\omega _{1}\right] _{\Omega
}\right) _{\Omega }.  \label{1.1}
\end{align}%
An integration by parts to right hand side then gives%
\begin{equation}
(\left. G_{0,1}(\omega _{1})\right\vert _{\Omega },\omega _{1})_{\Omega
}=\left\langle \frac{\partial \Delta \omega _{1}}{\partial n},\left[
R_{0}\omega _{1}\right] _{\Omega }\right\rangle _{\partial \Omega }-\left(
\nabla \Delta \omega _{1},\nabla \left[ R_{0}\omega _{1}\right] _{\Omega
}\right) _{\Omega }.  \label{1.2}
\end{equation}

\emph{To estimate the first term on the right hand side of (\ref{1.2}): }%
Using the trace estimate (\ref{trace2}), the regularity of the term $%
R_{0}\omega _{1}$ which is posted in (\ref{armed}) -- with therein $\delta
\equiv \frac{1}{2}$ -- and the Sobolev Trace Theorem, we have%
\begin{align*}
\left\vert \left\langle \frac{\partial \Delta \omega _{1}}{\partial n},\left[
R_{0}\omega _{1}\right] _{\Omega }\right\rangle _{\partial \Omega
}\right\vert & \leq \left\Vert \frac{\partial \Delta \omega _{1}}{\partial n}%
\right\Vert _{\partial \Omega }\left\Vert \left[ R_{0}\omega _{1}\right]
_{\Omega }\right\Vert _{\partial \Omega } \\
& \leq C\left\Vert \Delta \omega _{1}\right\Vert _{H^{1}(\Omega )}^{\frac{1}{%
2}}\left\Vert \Delta \omega _{1}\right\Vert _{H^{2}(\Omega )}^{\frac{1}{2}%
}\left\Vert R_{0}\omega _{1}\right\Vert _{H^{\frac{3}{2}}(\mathcal{O})} \\
& \leq C\left\Vert \omega _{1}\right\Vert _{H^{3}(\Omega )}^{\frac{1}{2}%
}\left\Vert \omega _{1}\right\Vert _{H^{4}(\Omega )}^{\frac{1}{2}}\left\Vert
\omega _{1}\right\Vert _{\Omega } \\
& \leq C\left\Vert \omega _{1}\right\Vert _{H^{2}(\Omega )}^{\frac{1}{4}%
}\left\Vert \omega _{1}\right\Vert _{H^{4}(\Omega )}^{\frac{1}{4}}\left\Vert
\omega _{1}\right\Vert _{H^{4}(\Omega )}^{\frac{1}{2}}\left\Vert \omega
_{1}\right\Vert _{\Omega } \\
& \leq C\left\Vert \omega _{1}\right\Vert _{H^{2}(\Omega )}^{\frac{1}{4}%
}\left\Vert [\omega _{1},\omega _{2},\mu ]\right\Vert _{D(\mathcal{A}_{0})}^{%
\frac{3}{4}}\left\Vert \omega _{1}\right\Vert _{\Omega }.
\end{align*}%
Using once more the resolvent relation (\ref{abstract}), we have for $%
\left\vert \alpha \right\vert >0$ sufficiently large,%
\begin{align}
\left\vert \left\langle \frac{\partial \Delta \omega _{1}}{\partial n},\left[
R_{\rho }\omega _{1}\right] _{\Omega }\right\rangle _{\partial \mathcal{O}%
}\right\vert & \leq C_{\vartheta }\left\vert \alpha \right\vert ^{\frac{3}{4}%
}\left\Vert \omega _{1}\right\Vert _{H^{2}(\Omega )}^{\frac{1}{4}}\left\Vert
[\omega _{1},\omega _{2},\mu ]-\frac{[\omega _{1}^{\ast },\omega _{2}^{\ast
},\mu ^{\ast }]}{\left\vert \alpha \right\vert e^{\pm i\vartheta }\sqrt{%
1+\tan ^{2}\vartheta }}\right\Vert _{\mathbf{H}_{0}}^{\frac{3}{4}}\left\Vert
\omega _{1}\right\Vert _{\Omega }  \notag \\
& \leq \epsilon \left\Vert \omega _{1}\right\Vert _{H^{2}(\Omega
)}^{2}+C_{\epsilon ,\vartheta }\left\vert \alpha \right\vert ^{\frac{6}{7}%
}\left\Vert [\omega _{1},\omega _{2},\mu ]-\frac{[\omega _{1}^{\ast },\omega
_{2}^{\ast },\mu ^{\ast }]}{\left\vert \alpha \right\vert e^{\pm i\vartheta }%
\sqrt{1+\tan ^{2}\vartheta }}\right\Vert _{\mathbf{H}_{0}}^{\frac{6}{7}%
}\left\Vert \omega _{1}\right\Vert _{\Omega }^{\frac{8}{7}}  \notag \\
& \leq C_{\epsilon ,\vartheta }\left\vert \alpha \right\vert ^{\frac{3}{2}%
}\left\Vert \omega _{1}\right\Vert _{\Omega }^{2}+2\epsilon \left\Vert
\lbrack \omega _{1},\omega _{2},\mu ]\right\Vert _{\mathbf{H}%
_{0}}^{2}+C_{\epsilon ,\vartheta }\left\Vert [\omega _{1}^{\ast },\omega
_{2}^{\ast },\mu ^{\ast }]\right\Vert _{\mathbf{H}_{0}}^{2},  \label{1.3}
\end{align}%
after again using Young's Inequality.

\bigskip

\emph{To estimate the second term on the right hand side of (\ref{1.2}):}
Using again the regularity of the term $R_{0}\omega _{1}$ which is posted in
(\ref{armed}) -- with therein $\delta \equiv \frac{1}{2}$, and the Sobolev
Trace Theorem, we have%
\begin{align*}
\left\vert \left( \nabla \Delta \omega _{1},\nabla \left[ R_{0}\omega _{1}%
\right] _{\Omega }\right) _{\Omega }\right\vert & \leq \left\Vert \nabla
\Delta \omega _{1}\right\Vert _{\Omega }\left\Vert \nabla \left[ R_{0}\omega
_{1}\right] _{\Omega }\right\Vert _{\Omega } \\
& \leq \left\Vert \nabla \Delta \omega _{1}\right\Vert _{\Omega }\left\Vert
\nabla \left[ R_{0}\omega _{1}\right] _{\Omega }\right\Vert _{H^{1}(\Omega )}
\\
& \leq C\left\Vert \omega _{1}\right\Vert _{H^{3}(\Omega )}\left\Vert
R_{0}\omega _{1}\right\Vert _{H^{\frac{3}{2}}(\mathcal{O})} \\
& \leq C\left\Vert \omega _{1}\right\Vert _{H^{2}(\Omega )}^{\frac{1}{2}%
}\left\Vert \omega _{1}\right\Vert _{H^{4}(\Omega )}^{\frac{1}{2}}\left\Vert
\omega _{1}\right\Vert _{\Omega } \\
& \leq C_{\vartheta }\left\Vert \omega _{1}\right\Vert _{H^{2}(\Omega )}^{%
\frac{1}{2}}\left\Vert [\omega _{1},\omega _{2},\mu ]\right\Vert _{D(%
\mathcal{A}_{0})}^{\frac{1}{2}}\left\Vert \omega _{1}\right\Vert _{\Omega }
\\
& \leq C_{\vartheta }\sqrt{\left\vert \alpha \right\vert }\left\Vert \omega
_{1}\right\Vert _{H^{2}(\Omega )}^{\frac{1}{2}}\left\Vert [\omega
_{1},\omega _{2},\mu ]-\frac{[\omega _{1}^{\ast },\omega _{2}^{\ast },\mu
^{\ast }]}{\left\vert \alpha \right\vert e^{\pm i\vartheta }\sqrt{1+\tan
^{2}\vartheta }}\right\Vert _{\mathbf{H}_{0}}^{\frac{1}{2}}\left\Vert \omega
_{1}\right\Vert _{\Omega },
\end{align*}%
after again using the resolvent relation (\ref{abstract}), and taking $%
\left\vert \alpha \right\vert >0$ sufficently large. Proceeding via Young's
Inequality, we have now%
\begin{align}
\left\vert \left( \nabla \Delta \omega _{1},\nabla \left[ R_{0}\omega _{1}%
\right] _{\Omega }\right) \right\vert & \leq \epsilon \left\Vert \omega
_{1}\right\Vert _{H^{2}(\Omega )}^{2}+C_{\epsilon ,\vartheta }\left\vert
\alpha \right\vert ^{\frac{2}{3}}\left\Vert [\omega _{1},\omega _{2},\mu ]-%
\frac{[\omega _{1}^{\ast },\omega _{2}^{\ast },\mu ^{\ast }]}{\left\vert
\alpha \right\vert e^{\pm i\vartheta }\sqrt{1+\tan ^{2}\vartheta }}%
\right\Vert _{\mathbf{H}_{0}}^{\frac{2}{3}}\left\Vert \omega _{1}\right\Vert
_{\Omega }^{\frac{4}{3}}  \notag \\
& \leq C_{\epsilon ,\vartheta }\left\vert \alpha \right\vert \left\Vert
\omega _{1}\right\Vert _{\Omega }^{2}+2\epsilon \left\Vert \lbrack \omega
_{1},\omega _{2},\mu ]\right\Vert _{\mathbf{H}_{0}}^{2}+C_{\epsilon
,\vartheta }\left\Vert [\omega _{1}^{\ast },\omega _{2}^{\ast },\mu ^{\ast
}]\right\Vert _{\mathbf{H}_{0}}^{2}.  \label{1.5}
\end{align}

\bigskip

Applying (\ref{1.3}) and (\ref{1.5}) to the right hand side of (\ref{1.2})
(and rescaling parameter $\epsilon >0$) now gives%
\begin{equation}
\left\vert (\left. G_{0,1}(\omega _{1})\right\vert _{\Omega },\omega
_{1})_{\Omega }\right\vert \leq C_{\epsilon ,\vartheta }\left\vert \alpha
\right\vert ^{\frac{3}{2}}\left\Vert \omega _{1}\right\Vert _{\Omega
}^{2}+2\epsilon \left\Vert \lbrack \omega _{1},\omega _{2},\mu ]\right\Vert
_{\mathbf{H}_{0}}^{2}+C_{\epsilon ,\vartheta }\left\Vert [\omega _{1}^{\ast
},\omega _{2}^{\ast },\mu ^{\ast }]\right\Vert _{\mathbf{H}_{0}}^{2}.
\label{1.7}
\end{equation}

\bigskip

\subsubsection{Analysis of the term $\left\vert \left( \left. G_{\protect%
\rho ,1}(\protect\omega _{1})\right\vert _{\Omega },\protect\omega %
_{1}\right) _{\Omega }\right\vert $ for $\protect\rho >0$}

Again from (\ref{G1}) and the expressions in (\ref{2.2}) and (\ref{2.4}), we
have

\begin{align}
(\left. G_{\rho ,1}(\omega _{1})\right\vert _{\Omega },\omega _{1})_{\Omega
}& =(\left. R_{\rho }P_{\rho }^{-1}\Delta ^{2}\omega _{1}\right\vert
_{\Omega },\omega _{1})_{\Omega }  \notag \\
& =\left( P_{\rho }^{-1}\Delta ^{2}\omega _{1},\left[ R_{\rho }\omega _{1}%
\right] _{\Omega }\right) _{\Omega }.  \label{1.8}
\end{align}%
At this point we reinvoke the positive definite, self-adjoint operator $%
\mathbf{\mathring{A}}:D(\mathring{A})\subset L^{2}(\Omega )\rightarrow
L^{2}(\Omega )$ in (\ref{ang}); with this operator in mind, we recall the
following characterizations (see \cite{grisvard}):%
\begin{equation}
D(\mathbf{\mathring{A}}^{\eta })\approx \left\{
\begin{array}{l}
\left\{ \varpi \in H^{4\eta }(\Omega ):\left. \varpi \right\vert _{\partial
\Omega }=0\right\} \text{, \ for }\frac{1}{8}<\eta <\frac{3}{8} \\
\\
\left\{ \varpi \in H^{4\eta }(\Omega ):\left. \varpi \right\vert _{\partial
\Omega }=\left. \frac{\partial \varpi }{\partial n}\right\vert _{\partial
\Omega }=0\right\} \text{, \ for }\frac{3}{8}<\eta \leq 1.%
\end{array}%
\right.  \label{char}
\end{equation}%
Proceeding from (\ref{1.8}) we have then,%
\begin{align}
(\left. (G_{\rho ,1}(\omega _{1})\right\vert _{\Omega },\omega _{1})_{\Omega
}& =\left( P_{\rho }^{-1}\mathbf{\mathring{A}}\omega _{1},\left[ R_{\rho
}\omega _{1}\right] _{\Omega }\right) _{\Omega }  \notag \\
& =\left( \mathbf{\mathring{A}}^{\frac{5}{8}+\epsilon }\omega _{1},\mathbf{%
\mathring{A}}^{\frac{3}{8}-\epsilon }P_{\rho }^{-1}\left[ R_{\rho }\omega
_{1}\right] _{\Omega }\right) _{\Omega }.  \label{1.9}
\end{align}

\bigskip

Using in part the fact that $\mathbf{\mathring{A}}^{\frac{3}{8}-\epsilon
}P_{\rho }^{-1}\left[ R_{\rho }\omega _{1}\right] _{\Omega }$ is in $%
L^{2}(\Omega ))$ continuously -- after using once more the estimate (\ref%
{armed}), with $\delta \equiv \frac{1}{2}$ -- a majorization of right hand
side then gives%
\begin{align*}
\left\vert (\left. (G_{\rho ,1}(\omega _{1})\right\vert _{\Omega },\omega
_{1})_{\Omega }\right\vert & \leq \left\Vert \mathbf{\mathring{A}}^{\frac{5}{%
8}+\epsilon }\omega _{1}\right\Vert _{\Omega }\left\Vert \mathbf{\mathring{A}%
}^{\frac{3}{8}-\epsilon }P_{\rho }^{-1}\left[ R_{\rho }\omega _{1}\right]
_{\Omega }\right\Vert _{\partial \Omega } \\
& \leq C\left\Vert \omega _{1}\right\Vert _{D(\mathbf{\mathring{A}}^{\frac{1%
}{2}})}^{\frac{1}{2}-4\epsilon }\left\Vert \omega _{1}\right\Vert _{D(%
\mathbf{\mathring{A}}^{\frac{3}{4}})}^{\frac{1}{2}+4\epsilon }\left\Vert
\omega _{1}\right\Vert _{\Omega } \\
& \leq C\left\Vert \omega _{1}\right\Vert _{D(\mathbf{\mathring{A}}^{\frac{1%
}{2}})}^{\frac{1}{2}-4\epsilon }\left\Vert [\omega _{1},\omega _{2},\mu
]\right\Vert _{D(\mathcal{A}_{\rho })}^{\frac{1}{2}+4\epsilon }\left\Vert
\omega _{1}\right\Vert _{H^{1}(\Omega )},
\end{align*}%
where in the last inequality we have recalled (\ref{S}) and (\ref{char}).
Using once more the resolvent relation (\ref{abstract}), we have for $%
\left\vert \alpha \right\vert >0$ sufficiently large,%
\begin{align}
& \left\vert (\left. (G_{\rho ,1}(\omega _{1})\right\vert _{\Omega },\omega
_{1})_{\Omega }\right\vert  \notag \\
& \leq C_{\vartheta }\left\vert \alpha \right\vert ^{\frac{1}{2}+4\epsilon
}\left\Vert \omega _{1}\right\Vert _{D(\mathbf{\mathring{A}}^{\frac{1}{2}%
})}^{\frac{1}{2}-4\epsilon }\left\Vert [\omega _{1},\omega _{2},\mu ]-\frac{%
[\omega _{1}^{\ast },\omega _{2}^{\ast },\mu ^{\ast }]}{\left\vert \alpha
\right\vert e^{\pm i\vartheta }\sqrt{1+\tan ^{2}\vartheta }}\right\Vert _{%
\mathbf{H}_{\rho }}^{\frac{1}{2}+4\epsilon }\left\Vert \omega
_{1}\right\Vert _{H^{1}(\Omega )}  \notag \\
& \leq \epsilon \left\Vert \omega _{1}\right\Vert _{D(\mathbf{\mathring{A}}^{%
\frac{1}{2}})}^{2}+C_{\epsilon ,\vartheta }\left\vert \alpha \right\vert ^{%
\frac{2+16\epsilon }{3+8\epsilon }}\left\Vert [\omega _{1},\omega _{2},\mu ]-%
\frac{[\omega _{1}^{\ast },\omega _{2}^{\ast },\mu ^{\ast }]}{\left\vert
\alpha \right\vert e^{\pm i\vartheta }\sqrt{1+\tan ^{2}\vartheta }}%
\right\Vert _{\mathbf{H}_{\rho }}^{\frac{2+16\epsilon }{3+8\epsilon }%
}\left\Vert \omega _{1}\right\Vert _{H^{1}(\Omega )}^{\frac{4}{3+8\epsilon }}
\notag \\
& \leq C_{\epsilon ,\vartheta }\left\vert \alpha \right\vert ^{1+8\epsilon
}\left\Vert \omega _{1}\right\Vert _{H^{1}(\Omega )}^{2}+2\epsilon
\left\Vert \lbrack \omega _{1},\omega _{2},\mu ]\right\Vert _{\mathbf{H}%
_{\rho }}^{2}+C_{\epsilon ,\vartheta }\left\Vert [\omega _{1}^{\ast },\omega
_{2}^{\ast },\mu ^{\ast }]\right\Vert _{\mathbf{H}_{\rho }}^{2},
\label{1.10}
\end{align}%
after using once more the characterization (\ref{char}).

\bigskip

Combining (\ref{1.7}) and (\ref{1.10}), we have then for all $\rho \geq 0$
and $\left\vert \alpha \right\vert >0$ sufficiently large,
\begin{equation}
\left\vert (\left. (G_{\rho ,1}(\omega _{1})\right\vert _{\Omega },\omega
_{1})_{\Omega }\right\vert \leq \left\vert (\left. G_{0,1}(\omega
_{1})\right\vert _{\Omega },\omega _{1})_{\Omega }\right\vert \leq
C_{\epsilon ,\vartheta }\left\vert \alpha \right\vert ^{\frac{3}{2}%
}\left\Vert \omega _{1}\right\Vert _{D(P_{\rho }^{\frac{1}{2}%
})}^{2}+\epsilon \left\Vert \lbrack \omega _{1},\omega _{2},\mu ]\right\Vert
_{\mathbf{H}_{\rho }}^{2}+C_{\epsilon ,\vartheta }\left\Vert [\omega
_{1}^{\ast },\omega _{2}^{\ast },\mu ^{\ast }]\right\Vert _{\mathbf{H}_{\rho
}}^{2}.  \label{G2.2}
\end{equation}

\subsection{The Proof Proper of Theorem \protect\ref{back}}

\bigskip

Applying the estimates (\ref{G22}) and (\ref{G2.2}), to the right hand side
of the expression (\ref{press}), and using the resolvent relation $\lambda
\omega _{1}=\omega _{2}+\omega _{1}^{\ast }$ yield the following lemma:

\begin{lemma}
For $\rho \geq 0$, the solution variables $[\omega _{1},\omega _{2},\mu ]$
of the resolvent equation obey the following estimate, for $\left\vert
\alpha \right\vert >0$ sufficiently large:%
\begin{align}
\left\vert \left\vert \left( \left. \pi _{0}\right\vert _{\Omega },\omega
_{1}\right) _{\Omega }\right\vert \right\vert & \leq C_{\epsilon, \vartheta
}\left\vert \alpha \right\vert ^{\frac{3}{2}+\epsilon }\left\Vert \omega
_{1}\right\Vert _{D(P_{\rho }^{\frac{1}{2}})}^{2}+\frac{\epsilon }{%
\left\vert \alpha \right\vert }\left\Vert \nabla \mu \right\Vert _{\mathcal{O%
}}^{2}  \notag \\
& \,\,\,\,+\left( \epsilon +\frac{C_{\epsilon, \vartheta }}{\left\vert
\alpha \right\vert ^{\frac{1}{2}-\epsilon }}\right) \left( \left\Vert
[\omega _{1},\omega _{2},\mu ]\right\Vert _{\mathbf{H}_{\rho
}}^{2}+C_{\epsilon, \vartheta }\left\Vert [\omega _{1}^{\ast },\omega
_{2}^{\ast },\mu ^{\ast }]\right\Vert _{\mathbf{H}_{\rho }}^{2}\right)
\notag \\
&  \notag \\
& \leq \frac{\epsilon }{\left\vert \alpha \right\vert }\left\Vert \nabla \mu
\right\Vert _{\mathcal{O}}^{2}+\left( \epsilon +\frac{C_{\epsilon ,\vartheta
}}{\left\vert \alpha \right\vert ^{\frac{1}{2}-\epsilon }}\right) \left\Vert
[\omega _{1},\omega _{2},\mu ]\right\Vert _{\mathbf{H}_{\rho
}}^{2}+C_{\epsilon ,\vartheta }\left\Vert [\omega _{1}^{\ast },\omega
_{2}^{\ast },\mu ^{\ast }]\right\Vert _{\mathbf{H}_{\rho }}^{2}.
\label{pressure_cross}
\end{align}
\end{lemma}

\bigskip

In completing the proof of Theorem \ref{back}, we bear in mind that Criteria
1 and 2 are imposed upon complex parameter $\lambda =\alpha +i\beta $.

\bigskip

\emph{Step 1.} We apply the estimate (\ref{pressure_cross}) to the right
hand side of (\ref{p1.2}), so as to have%
\begin{align}
\left\vert \beta \right\vert \left\Vert \mu \right\Vert _{\mathcal{O}}^{2}&
=\left\vert -\text{Im}\left( \left. \pi _{0}\right\vert _{\Omega },\lambda
\omega _{1}-\omega _{1}^{\ast }\right) _{\Omega }+\text{Im}\left( \mu ^{\ast
},\mu \right) _{\mathcal{O}}\right\vert  \notag \\
& \leq C_{\vartheta }\left\vert \alpha \right\vert \left\vert \text{Im}%
\left( \left. \pi _{0}\right\vert _{\Omega },\omega _{1}\right) _{\Omega
}\right\vert +\left\vert \text{Im}\left( \left. \pi _{0}\right\vert _{\Omega
},\omega _{1}^{\ast }\right) _{\Omega }\right\vert +\left\vert \text{Im}%
\left( \mu ^{\ast },\mu \right) \right\vert _{\mathcal{O}}  \notag \\
& \leq \left\vert \text{Im}\left( \left. \pi _{0}\right\vert _{\Omega
},\omega _{1}^{\ast }\right) _{\Omega }\right\vert +\epsilon C_{\vartheta
}^{\ast }\left\Vert \nabla \mu \right\Vert _{\mathcal{O}}^{2}  \notag \\
& \,\,\,\,\,+C_{\vartheta }^{\ast }\left\vert \alpha \right\vert \left(
\epsilon +\frac{C_{\epsilon }}{\left\vert \alpha \right\vert ^{\frac{1}{2}%
-\epsilon }}\right) \left\Vert [\omega _{1},\omega _{2},\mu ]\right\Vert _{%
\mathbf{H}_{\rho }}^{2}+C_{\epsilon ,\vartheta }\left\vert \alpha
\right\vert \left\Vert [\omega _{1}^{\ast },\omega _{2}^{\ast },\mu ^{\ast
}]\right\Vert _{\mathbf{H}_{\rho }}^{2},  \label{s2}
\end{align}%
where above, positive constant $C_{\vartheta }^{\ast }$ is independent of
parameter $\epsilon >0.$

\medskip

Now as for the first term on the right hand side of (\ref{s2}): Since the
datum $\omega _{1}^{\ast }$ satisfies the compatibility condition $%
\int_{\Omega }\omega _{1}^{\ast }d\Omega =0$, then there is a function $%
\varphi (\omega _{1}^{\ast })\in \mathbf{H}^{1}(\mathcal{O})$ which solves%
\begin{align}
\text{div}(\varphi )& =0\text{ \ in }\mathcal{O}\text{;} & &  \notag
\label{div} \\
\varphi & =0\text{ \ in }S; & &  \notag \\
\varphi & =[0,0,\omega _{1}^{\ast }]\text{ \ in }\Omega , & &
\end{align}%
with the estimate%
\begin{equation}
\left\Vert \nabla \varphi \right\Vert _{\mathcal{O}}\leq C\left\Vert \omega
_{1}^{\ast }\right\Vert _{H^{\frac{1}{2}}(\Omega )}  \label{div2}
\end{equation}%
(see e.g., p. 9 of \cite{galdi}). With this solution variable $\varphi
(\omega _{1}^{\ast })$ in hand, and by virtue of the geometry in play, we
then have%
\begin{eqnarray*}
\left( \left. \pi _{0}\right\vert _{\Omega },\omega _{1}^{\ast }\right)
_{\Omega } &=&-\left( \frac{\partial \mu }{\partial \nu },\varphi \right)
_{\partial \mathcal{O}}+\left( \pi _{0}\nu ,\varphi \right) _{\partial
\mathcal{O}} \\
&&-\left( \nabla \mu ,\nabla \varphi \right) _{\mathcal{O}}-\left( \Delta
\mu ,\varphi \right) _{\mathcal{O}}+\left( \nabla \pi _{0},\varphi \right)
_{\partial \mathcal{O}}+0 \\
&=&-\left( \nabla \mu ,\nabla \varphi \right) _{\mathcal{O}}-\lambda \left(
\mu ,\varphi \right) _{\mathcal{O}}+\left( \mu ^{\ast },\varphi \right) _{%
\mathcal{O}},
\end{eqnarray*}%
after using the fluid equation in (\ref{f}). We have then upon majorizing,
with the use of the estimate (\ref{div2}), and for large $\left\vert \alpha
\right\vert >0$%
\begin{equation}
\left\vert \left( \left. \pi _{0}\right\vert _{\Omega },\omega _{1}^{\ast
}\right) _{\Omega }\right\vert \leq \epsilon \left( \left\Vert \nabla \mu
\right\Vert _{\mathcal{O}}^{2}+\left\vert \alpha \right\vert \left\Vert \mu
\right\Vert _{\mathcal{O}}^{2}\right) +\left\vert \alpha \right\vert
C_{\epsilon ,\vartheta }\left\Vert [\omega _{1}^{\ast },\omega _{2}^{\ast
},\mu ^{\ast }]\right\Vert _{\mathbf{H}_{\rho }}^{2}.  \label{s2.5}
\end{equation}%
Applying this estimate to the right hand side of (\ref{s2}) now yields
(after a rescaling of $\epsilon >0$)%
\begin{align}
\left\vert \beta \right\vert \left\Vert \mu \right\Vert _{\mathcal{O}}^{2}&
\leq \epsilon C_{\vartheta }^{\ast }\left\Vert \nabla \mu \right\Vert _{%
\mathcal{O}}^{2}+C_{\vartheta }^{\ast }\left\vert \alpha \right\vert \left(
\epsilon +\frac{C_{\epsilon }}{\left\vert \alpha \right\vert ^{\frac{1}{2}%
-\epsilon }}\right) \left\Vert [\omega _{1},\omega _{2},\mu ]\right\Vert _{%
\mathbf{H}_{\rho }}^{2}  \notag \\
& +C_{\epsilon ,\vartheta }\left\vert \alpha \right\vert \left\Vert [\omega
_{1}^{\ast },\omega _{2}^{\ast },\mu ^{\ast }]\right\Vert _{\mathbf{H}_{\rho
}}^{2},  \label{s3}
\end{align}%
where again, positive constant $C_{\vartheta }^{\ast }$ is independent of
parameter $\epsilon >0.$

\bigskip

\emph{Step 2:} We invoke the relation (\ref{1.22}):%
\begin{equation*}
\left\Vert \nabla \mu \right\Vert _{\mathcal{O}}^{2}=\left\vert \alpha
\right\vert \left\Vert \mu \right\Vert _{\mathcal{O}}^{2}-\text{Re}\left(
\left. \pi _{0}\right\vert _{\Omega },\lambda \omega _{1}-\omega _{1}^{\ast
}\right) _{\Omega }+\text{Re}\left( \mu ^{\ast },\mu \right) _{\mathcal{O}}.
\end{equation*}%
Applying the estimates (\ref{s3}), (\ref{pressure_cross}), and (\ref{s2.5})
to right hand side now gives%
\begin{align}
\left\Vert \nabla \mu \right\Vert _{\mathcal{O}}^{2}& \leq \epsilon
C_{\vartheta }^{\ast }\left\Vert \nabla \mu \right\Vert _{\mathcal{O}%
}^{2}+C_{\vartheta }^{\ast }\left\vert \alpha \right\vert \left( \epsilon +%
\frac{C_{\epsilon }}{\left\vert \alpha \right\vert ^{\frac{1}{2}-\epsilon }}%
\right) \left\Vert [\omega _{1},\omega _{2},\mu ]\right\Vert _{\mathbf{H}%
_{\rho }}^{2}  \notag \\
& +C_{\epsilon ,\vartheta }\left\vert \alpha \right\vert \left\Vert [\omega
_{1}^{\ast },\omega _{2}^{\ast },\mu ^{\ast }]\right\Vert _{\mathbf{H}_{\rho
}}^{2}.  \label{s4}
\end{align}

\bigskip

\emph{Step 3:} We apply the estimate (\ref{s4}) to the right hand side of (%
\ref{p1.3}). This gives for large $\left\vert \alpha \right\vert >0$%
\begin{align}
\left\vert \alpha \right\vert \left\Vert \left[
\begin{array}{c}
\omega _{1} \\
\omega _{2} \\
\mu%
\end{array}%
\right] \right\Vert _{\mathbf{H}_{\rho }}^{2}& =\left\vert \left\Vert \nabla
\mu \right\Vert _{\mathcal{O}}^{2}+\text{Re}\left( \left[
\begin{array}{c}
\omega _{1}^{\ast } \\
\omega _{2}^{\ast } \\
\mu ^{\ast }%
\end{array}%
\right] ,\left[
\begin{array}{c}
\omega _{1} \\
\omega _{2} \\
\mu%
\end{array}%
\right] \right) _{\mathbf{H}_{\rho }}\right\vert  \notag \\
& \leq \epsilon C_{\vartheta }^{\ast }\left\Vert \nabla \mu \right\Vert _{%
\mathcal{O}}^{2}+C_{\vartheta }^{\ast }\left\vert \alpha \right\vert \left(
\epsilon +\frac{C_{\epsilon }}{\left\vert \alpha \right\vert ^{\frac{1}{2}%
-\epsilon }}\right) \left\Vert [\omega _{1},\omega _{2},\mu ]\right\Vert _{%
\mathbf{H}_{\rho }}^{2}  \notag \\
& +C_{\epsilon ,\vartheta }\left\vert \alpha \right\vert \left\Vert [\omega
_{1}^{\ast },\omega _{2}^{\ast },\mu ^{\ast }]\right\Vert _{\mathbf{H}_{\rho
}}^{2}.  \label{s5}
\end{align}

\bigskip

\emph{Step 4:} Taking $\epsilon >0$ sufficently small in (\ref{s5}) (with
again positive constant $C_{\vartheta }^{\ast }$ being independent of
parameter $\epsilon >0$), we have%
\begin{eqnarray*}
&&\left( 1-\epsilon C_{\vartheta }^{\ast }\right) \left\vert \alpha
\right\vert \left\Vert \left[
\begin{array}{c}
\omega _{1} \\
\omega _{2} \\
\mu%
\end{array}%
\right] \right\Vert _{\mathbf{H}_{\rho }}^{2} \\
&\leq &\epsilon C_{\vartheta }^{\ast }\left\Vert \nabla \mu \right\Vert _{%
\mathcal{O}}^{2}+C_{\epsilon ,\vartheta }\frac{\left\vert \alpha \right\vert
}{\left\vert \alpha \right\vert ^{\frac{1}{2}-\epsilon }}\left\Vert \left[
\begin{array}{c}
\omega _{1} \\
\omega _{2} \\
\mu%
\end{array}%
\right] \right\Vert _{\mathbf{H}_{\rho }}^{2}+C_{\epsilon ,\vartheta
}\left\vert \alpha \right\vert \left\Vert \left[
\begin{array}{c}
\omega _{1}^{\ast } \\
\omega _{2}^{\ast } \\
\mu ^{\ast }%
\end{array}%
\right] \right\Vert _{\mathbf{H}_{\rho }}^{2},
\end{eqnarray*}%
and so
\begin{align}
& \left\vert \alpha \right\vert \left\Vert \left[
\begin{array}{c}
\omega _{1} \\
\omega _{2} \\
\mu%
\end{array}%
\right] \right\Vert _{\mathbf{H}_{\rho }}^{2}  \notag \\
& \leq \frac{\epsilon C_{\vartheta }^{\ast }}{1-\epsilon C_{\vartheta
}^{\ast }}\left\Vert \nabla \mu \right\Vert _{\mathcal{O}}^{2}+C_{\epsilon
,\vartheta }\frac{\left\vert \alpha \right\vert }{\left\vert \alpha
\right\vert ^{\frac{1}{2}-\epsilon }}\left\Vert \left[
\begin{array}{c}
\omega _{1} \\
\omega _{2} \\
\mu%
\end{array}%
\right] \right\Vert _{\mathbf{H}_{\rho }}^{2}+C_{\epsilon ,\vartheta
}\left\vert \alpha \right\vert \left\Vert \left[
\begin{array}{c}
\omega _{1}^{\ast } \\
\omega _{2}^{\ast } \\
\mu ^{\ast }%
\end{array}%
\right] \right\Vert _{\mathbf{H}_{\rho }}^{2}.  \label{s6}
\end{align}

\bigskip

\emph{Step 5:} We return to the estimate (\ref{s4}). Applying (\ref{s6})
thereto gives for $\epsilon >0$ sufficently small,%
\begin{align}
\left\Vert \nabla \mu \right\Vert _{\mathcal{O}}^{2}& \leq \epsilon
C_{\vartheta }^{\ast }\left\Vert \nabla \mu \right\Vert _{\mathcal{O}}^{2}+%
\frac{C_{\epsilon ,\vartheta }\left\vert \alpha \right\vert }{\left\vert
\alpha \right\vert ^{\frac{1}{2}-\epsilon }}\left\Vert [\omega _{1},\omega
_{2},\mu ]\right\Vert _{\mathbf{H}_{\rho }}^{2}  \notag \\
& +\,\,C_{\epsilon ,\vartheta }\left\vert \alpha \right\vert \left\Vert
[\omega _{1}^{\ast },\omega _{2}^{\ast },\mu ^{\ast }]\right\Vert _{\mathbf{H%
}_{\rho }}^{2},  \label{s7}
\end{align}%
where again (relabeled) positive constant $C_{\vartheta }^{\ast }$ does not
depend upon small $\epsilon >0$. Further specifying $\epsilon >0$ to be
small enough, we have now%
\begin{equation}
\left( 1-\epsilon C_{\vartheta }^{\ast }\right) \left\Vert \nabla \mu
\right\Vert _{\mathcal{O}}^{2}\leq \frac{C_{\epsilon ,\vartheta }\left\vert
\alpha \right\vert }{\left\vert \alpha \right\vert ^{\frac{1}{2}-\epsilon }}%
\left\Vert [\omega _{1},\omega _{2},\mu ]\right\Vert _{\mathbf{H}_{\rho
}}^{2}+C_{\epsilon ,\vartheta }\left\vert \alpha \right\vert \left\Vert
[\omega _{1}^{\ast },\omega _{2}^{\ast },\mu ^{\ast }]\right\Vert _{\mathbf{H%
}_{\rho }}^{2},
\end{equation}%
whence we obtain%
\begin{equation}
\left\Vert \nabla \mu \right\Vert _{\mathcal{O}}^{2}\leq \frac{C_{\epsilon
,\vartheta }\left\vert \alpha \right\vert }{\left\vert \alpha \right\vert ^{%
\frac{1}{2}-\epsilon }}\left\Vert [\omega _{1},\omega _{2},\mu ]\right\Vert
_{\mathbf{H}_{\rho }}^{2}+C_{\epsilon ,\vartheta }\left\vert \alpha
\right\vert \left\Vert [\omega _{1}^{\ast },\omega _{2}^{\ast },\mu ^{\ast
}]\right\Vert _{\mathbf{H}_{\rho }}^{2}.  \label{s8}
\end{equation}

\bigskip

\emph{Step 6:} We finish the proof by applying the estimate (\ref{s8}) to
the right hand side of (\ref{s6}). This gives%
\begin{equation}
\left\vert \alpha \right\vert \left\Vert \left[
\begin{array}{c}
\omega _{1} \\
\omega _{2} \\
\mu%
\end{array}%
\right] \right\Vert _{\mathbf{H}_{\rho }}^{2}\leq C_{\epsilon ,\vartheta }%
\frac{\left\vert \alpha \right\vert }{\left\vert \alpha \right\vert ^{\frac{1%
}{2}-\epsilon }}\left\Vert \left[
\begin{array}{c}
\omega _{1} \\
\omega _{2} \\
\mu%
\end{array}%
\right] \right\Vert _{\mathbf{H}_{\rho }}^{2}+C_{\epsilon ,\vartheta
}\left\vert \alpha \right\vert \left\Vert \left[
\begin{array}{c}
\omega _{1}^{\ast } \\
\omega _{2}^{\ast } \\
\mu ^{\ast }%
\end{array}%
\right] \right\Vert _{\mathbf{H}_{\rho }}^{2},  \label{s9}
\end{equation}%
for fixed $\epsilon >0$, small enough. Taking now $\left\vert \alpha
\right\vert $ so large that $1-\frac{C_{\epsilon ,\vartheta }}{\left\vert
\alpha \right\vert ^{\frac{1}{2}-\epsilon }}>\frac{1}{2}$; i.e.,%
\begin{equation*}
\left\vert \alpha \right\vert >\left( 2C_{\epsilon ,\vartheta }\right) ^{%
\frac{2}{1-2\epsilon }},
\end{equation*}%
we have finally%
\begin{equation}
\frac{\left\vert \alpha \right\vert }{2}\left\Vert \left[
\begin{array}{c}
\omega _{1} \\
\omega _{2} \\
\mu%
\end{array}%
\right] \right\Vert _{\mathbf{H}_{\rho }}^{2}\leq C_{\epsilon ,\vartheta
}\left\vert \alpha \right\vert \left\Vert \left[
\begin{array}{c}
\omega _{1}^{\ast } \\
\omega _{2}^{\ast } \\
\mu ^{\ast }%
\end{array}%
\right] \right\Vert _{\mathbf{H}_{\rho }}^{2},  \label{final}
\end{equation}%
which gives the uniform bound (\ref{bound}). This completes the proof of
Theorem \ref{back}.

\end{document}